\documentclass[12pt]{article}
\usepackage{amssymb,amsmath}
\begin{document}
\begin{center}
{\Large\bf Interior capacities of~condensers\\ in~locally compact spaces}\\
\bigskip
{\large\it Natalia Zorii}\\
\end{center}
\bigskip

{\small{\bf Abstract.} The study is motivated by the known fact
that, in the noncompact case, the main minimum-problem of the theory
of interior capacities of condensers in a locally compact space is
in general unsolvable, and this occurs even under very natural
assumptions (e.\,g., for the Newton, Green, or Riesz kernels
in~$\mathbb R^n$ and closed condensers). Therefore it was
particularly interesting to find statements of variational problems
dual to the main min\-imum-prob\-lem (and hence providing some new
equivalent definitions of the capacity), but always solvable
(e.\,g., even for nonclosed condensers). For all positive definite
kernels satisfying B.~Fuglede's condition of consistency between the
strong and vague topologies, problems with the desired properties
are posed and solved. Their solutions provide a natural
generalization of the well-known notion of interior capacitary
distributions associated with a set. We give a description of those
solutions, establish statements on their uniqueness and continuity,
and point out their characteristic properties. A condenser is
treated as a finite collection of arbitrary sets with sing~$+1$
or~$-1$ prescribed, such that the closures of opposite-signed sets
are mutually disjoint.\bigskip

\bigskip

}\bigskip

{\bf\large 1. Introduction}\medskip

The present work is devoted to further development of the theory of
interior capacities of arbitrary (noncompact or even nonclosed)
condensers in a locally compact space, started by the author
in~\cite{Z3,Z4}. For a background of that theory in the compact
case, see the study by~M.~Ohtsuka~\cite{O}.

The reader is expected to be familiar with the principal notions and
results of the theory of measures and integration on a locally
compact space; its exposition can be found in~\cite{B2,E2} (see
also~\cite{F1,Z3} for a brief survey).

The theory of interior capacities of condensers provides a natural
extension of the well-known theory of interior capacities of sets,
developed by H.~Cartan~\cite{Car} and Vall\'{e}e-Poussin~\cite{VP}
for classical kernels in~$\mathbb R^n$ and later on generalized by
B.~Fuglede~\cite{F1} for general kernels in a locally compact
space~$\mathbf X$. However, those two theories
--- for sets and, on the other hand, condensers~--- are drastically
different. To illustrate this, it is enough to note that, in the
noncompact case, the main mi\-ni\-mum-pro\-blem of the theory of
interior capacities of condensers is in general {\it unsolvable\/},
and this phenomenon occurs even under very natural assumptions
(e.\,g., for the Newton, Green, or Riesz kernels in~$\mathbb R^n$
and closed condensers); compare with~\cite{Car, F1}. Necessary and
sufficient conditions for the problem to be solvable were given
in~\cite{Z4,Z6}; see~Sec.~4.5 below for a brief survey.

Therefore it was particularly interesting to find statements of
variational problems dual to the main minimum-problem of the theory
of interior capacities of condensers, but in contrast to the last
one, {\it always\/} solvable~--- e.\,g., even for nonclosed
condensers. (When speaking on duality of variational problems, we
always mean their extremal values to be equal.)

In all that follows, $\mathbf X$ denotes a locally compact Hausdorff
space, and $\mathfrak M=\mathfrak M(\mathbf X)$ the linear space of
all real-valued Radon measures~$\nu$ on~$\mathbf X$ equipped with
the {\it vague\/} topology, i.\,e., the topology of pointwise
convergence on the class $\mathbf C_0(\mathbf X)$ of all real-valued
continuous functions on~$\mathbf X$ with compact support.

A {\it kernel\/} $\kappa$ on $\mathbf X$ is meant to be a lower
semicontinuous function $\kappa:\mathbf X\times\mathbf
X\to(-\infty,\infty]$. To avoid some difficulties when defining
energies and potentials, we follow~\cite{F1} in assuming that
$\kappa\geqslant0$ unless the space~$\mathbf X$ is compact.

The {\it energy\/} and the {\it potential\/} of a measure
$\nu\in\mathfrak M$ with respect to a kernel~$\kappa$ are defined by
$$\kappa(\nu,\nu):=\int\kappa(x,y)\,d(\nu\otimes\nu)(x,y)$$
and $$\kappa(x,\nu):=\int\kappa(x,y)\,d\nu(y),\quad x\in\mathbf X,$$
respectively, provided the corresponding integral above is well
defined (as a finite number or $\pm\infty$). Let $\mathcal E$ denote
the set of all $\nu\in\mathfrak M$ with
$-\infty<\kappa(\nu,\nu)<\infty$.

In the present study we shall be concerned with minimal energy
problems over certain subclasses of~$\mathcal E$, properly chosen.
For all positive definite kernels satisfying B.~Fuglede's condition
of consistency between the strong and vague topologies on~$\mathcal
E$ (see~Sec.~2 below), those variational problems are shown to be
{\it dual\/} to the main minimum-pro\-blem of the theory of interior
capacities of con\-den\-sers (and hence providing some new {\it
equivalent\/} definitions of the capacity), but {\it always
solvable\/}. See~Theorems~2\,--\,4 and Corollaries~6,~8.

Their solutions provide a natural generalization of the notion of
interior capacitary distributions associated with a set, introduced
in~\cite{F1}. We shall give a description of those solutions,
establish statements on their uniqueness and continuity, and point
out their characteristic properties (see~Sec.~6\,--\,9).

Condensers and their capacities are treated in a fairly general
sense; see~Sec.~3 and~4 below for the corresponding
definitions.\bigskip

{\bf\large 2. Preliminaries: topologies, consistent and perfect
kernels}\nopagebreak\medskip

Recall that a measure $\nu\geqslant0$ is said to be {\it
concentrated\/} on~$E$, where $E$ is a subset of~$\mathbf X$, if the
complement $\complement E:=\mathbf X\setminus E$ is locally
$\nu$-negligible; or, equivalently, if $E$ is $\nu$-measurable and
$\nu=\nu_E$, where $\nu_E$ denotes the trace of~$\nu$ upon~$E$.

We denote by $\mathfrak M^+(E)$ the convex cone of all nonnegative
measures concentrated on~$E$, and write $\mathcal E^+(E):=\mathfrak
M^+(E)\cap\mathcal E$. To shorten notation, write $\mathcal
E^+:=\mathcal E^+(\mathbf X)$.

From now on, the kernel under consideration is always assumed to be
{\it positive definite}, which means that it is symmetric (i.\,e.,
$\kappa(x,y)=\kappa(y,x)$ for all $x,\,y\in\mathbf X$) and the
energy $\kappa(\nu,\nu)$, $\nu\in\mathfrak M$, is nonnegative
whenever defined. Then $\mathcal E$ is known to be a pre-Hil\-bert
space with the scalar product
$$
\kappa(\nu_1,\nu_2):=\int\kappa(x,y)\,d(\nu_1\otimes\nu_2)(x,y)
$$
and the seminorm $\|\nu\|:=\sqrt{\kappa(\nu,\nu)}$; see~\cite{F1}. A
(positive definite) kernel is called {\it strictly positive
definite\/} if the seminorm $\|\cdot\|$ is a norm.

A measure $\nu\in\mathcal E$ is said to be {\it equivalent in\/}
$\mathcal E$ to a given $\nu_0\in\mathcal E$ if $\|\nu-\nu_0\|=0$;
the equivalence class, consisting of all those~$\nu$, will be
denoted by~$\left[\nu_0\right]_\mathcal E$.

In addition to the {\it strong\/} topology on~$\mathcal E$,
determined by the above seminorm~$\|\cdot\|$, it is often useful to
consider the {\it weak\/} topology on~$\mathcal E$, defined by means
of the seminorms $\nu\mapsto|\kappa(\nu,\mu)|$, $\mu\in\mathcal E$
(see~\cite{F1}). The Cauchy-Schwarz inequality
$$
|\kappa(\nu,\mu)|\leqslant\|\nu\|\,\|\mu\|,\quad\nu,\,\mu\in\mathcal
E,
$$
implies immediately that the strong topology on $\mathcal E$ is
finer than the weak one.

In \cite{F1}, B.~Fuglede introduced the following two properties of
{\it consistency\/} between the induced strong, weak, and vague
topologies on~$\mathcal E^+$:\smallskip

($C$) \ {\it Every strong Cauchy net in $\mathcal E^+$ converges
strongly to every its vague cluster point;}\smallskip

($CW$) \ {\it Every strongly bounded and vaguely convergent net in
$\mathcal E^+$ converges weakly to the vague limit;}\smallskip

\noindent in \cite{F2}, the properties ($C$) and ($CW$) were shown
to be {\it equivalent\/}.\medskip

{\bf Definition 1.} Following B.~Fuglede, we call a kernel~$\kappa$
{\it consistent} if it satisfies either of the properties~($C$)
and~($CW$), and {\it perfect\/} if, in addition, it is strictly
positive definite.\medskip

{\bf Remark 1.} One has to consider {\it nets\/} or {\it filters\/}
in~$\mathfrak M$ instead of sequences, because the vague topology in
general does not satisfy the first axiom of countability. We follow
Moore's and Smith's theory of convergence, based on the concept of
nets (see~\cite{MS}; cf.~also~\cite[Ch.~0]{E2} and
\cite[Ch.~2]{K}).\medskip

{\bf Theorem 1}~\cite{F1}. {\it A kernel $\kappa$ is perfect if and
only if $\mathcal E^+$ is strongly complete and the strong topology
on~$\mathcal E^+$ is finer than the vague one.}\medskip

{\bf Examples.} In $\mathbb R^n$, $n\geqslant 3$, the Newton kernel
$|x-y|^{2-n}$ is perfect~\cite{Car}. So is the Riesz kernel
$|x-y|^{\alpha-n}$, $0<\alpha<n$, in~$\mathbb R^n$, $n\geqslant2$
(see~\cite{D1, D2}). Furthermore, if $D$ is an open set in~$\mathbb
R^n$, $n\geqslant 2$, and its generalized Green function~$g_D$
exists (see, e.\,g.,~\cite[Th.~5.24]{HK}), then the Green
kernel~$g_D$ is perfect as well~\cite{E1}.\medskip

{\bf Remark 2.} As is seen from Theorem~1, the concept of consistent
or perfect kernels is an efficient tool in minimal energy problems
over classes of {\it nonnegative\/} measures with finite energy.
Indeed, the theory of capacities of {\it sets\/} has been developed
in~\cite{F1} exactly for those kernels. We shall show below that
this concept is still efficient in minimal energy problems over
classes of {\it signed\/} measures associated with a {\it
condenser}. This is guaranteed by the theorem on the strong
completeness of proper subspaces of~$\mathcal E$, to be stated in
Sec.~10 below.\bigskip

{\bf\large 3. Condensers. Measures associated with a condenser;
their energies and potentials}\nopagebreak\medskip

{\bf 3.1.} Fix natural numbers $m$ and $m+p$, where $p\geqslant0$,
and write
$$
I:=\bigl\{1,\ldots,m+p\bigr\},\qquad
I^+:=\bigl\{1,\ldots,m\bigr\},\qquad I^-:=I\setminus I^+.
$$

{\bf Definition 2.} An ordered collection $\mathcal A=(A_i)_{i\in
I}$ of nonempty sets $A_i\subset\mathbf X$, $i\in I$, is called an
$(m,p)$-{\it condenser\/} in~$\mathbf X$ (or simply a {\it
condenser\/}) if
\begin{equation}
\overline{A_i}\cap\overline{A_j}=\varnothing\quad\mbox{for all}\quad
i\in I^+,\quad j\in I^-. \label{non}
\end{equation}

The sets $A_i$, $i\in I^+$, and $A_j$, $j\in I^-$, are said to be
the {\it positive\/} and, respectively, {\it negative plates\/} of
an $(m,p)$-condenser $\mathcal A=(A_i)_{i\in I}$. Note that any two
equal-sign\-ed plates of a condenser are allowed to intersect each
other.

Let $\mathfrak C_{m,p}=\mathfrak C_{m,p}(\mathbf X)$ denote the
collection of all $(m,p)$-condensers in~$\mathbf X$. A condenser
$\mathcal A\in\mathfrak C_{m,p}$ is called {\it closed\/} or {\it
compact\/} if all the plates~$A_i$, $i\in I$, are closed or,
respectively, compact. Similarly, we shall call it {\it universally
measurable\/} if all the plates are universally measurable~--- that
is, measurable with respect to every nonnegative Radon measure.

Given $\mathcal A=(A_i)_{i\in I}$, write $\overline{\mathcal
A}:=(\,\overline{A_i}\,)_{i\in I}$. Then, due to~(\ref{non}),
$\overline{\mathcal A}$ is a (closed) $(m,p)$-condenser. In the
sequel, also the following notation will be required:
$$
A:=\bigcup_{i\in I}\,A_i,\qquad A^+:=\bigcup_{i\in I^+}\,A_i,\qquad
A^-:=\bigcup_{i\in I^-}\,A_i.
$$

{\bf 3.2.} With the preceding notation, write
$$
\alpha_i:=\left\{
\begin{array}{rll} +1 & \mbox{if} & i\in I^+,\\ -1 & \mbox{if} & i\in
I^-.\\ \end{array} \right.
$$
Given $\mathcal A\in\mathfrak C_{m,p}$, let $\mathfrak M(\mathcal
A)$ consist of all {\it linear combinations\/} of the form
$$
\mu=\sum_{i\in I}\,\alpha_i\mu^i,\mbox{ \ where \ }\mu^i\in\mathfrak
M^+(A_i)\mbox{ \ for all \ } i\in I.
$$
Any two $\mu_1$ and $\mu_2$ in $\mathfrak M(\mathcal A)$,
$$\mu_1=\sum_{i\in
I}\,\alpha_i\mu_1^i\quad\mbox{and}\quad\mu_2=\sum_{i\in
I}\,\alpha_i\mu_2^i,
$$
are regarded to be {\it identical\/} ($\mu_1\equiv\mu_2$) if and
only if $\mu_1^i=\mu_2^i$ for all $i\in I$.

Note that, under the relation of identity in $\mathfrak M(\mathcal
A)$ thus defined, the following correspondence is one-to-one:
$$
\mathfrak M(\mathcal A)\ni\mu\mapsto(\mu^i)_{i\in I}\in \prod_{i\in
I}\mathfrak M^+(A_i).
$$
We shall call $\mu\in\mathfrak M(\mathcal A)$ a {\it measure
associated with a condenser\/}~$\mathcal A$, and the
measure~$\mu^i$, $i\in I$, its $i$-{\it coordinate}.

For any $\mu_1,\,\mu_2\in\mathfrak M(\mathcal A)$ and
$q_1,\,q_2\in\mathbb R_+$, define $q_1\mu_1+q_2\mu_2$ to be an
element from $\mathfrak M(\mathcal A)$ uniquely determined by the
relations
$$
(q_1\mu_1+q_2\mu_2)^i:=q_1\mu_1^i+q_2\mu_2^i,\quad i\in I.
$$
Then the set $\mathfrak M(\mathcal A)$ becomes {\it
convex\/}.\medskip

{\bf 3.3.} Given $\mu\in\mathfrak M(\mathcal A)$, denote by~$R\mu$
the Radon measure uniquely determined by either of the two
equivalent relations
$$
R\mu(\varphi)=\sum_{i\in I}\,\alpha_i\mu^i(\varphi)\quad\mbox{for
all \ }\varphi\in\mathbf C_0(\mathbf X),
$$
\begin{equation}
(R\mu)^+=\sum_{i\in I^+}\,\mu^i,\qquad(R\mu)^-=\sum_{i\in
I^-}\,\mu^i.\label{R}
\end{equation}

Of course, the mapping~$R:\mathfrak M(\mathcal A)\to\mathfrak M$
thus defined is in general non-injective, i.\,e., one can choose
$\mu'\in\mathfrak M(\mathcal A)$ so that $\mu'\not\equiv\mu$, while
$R\mu'=R\mu$. (It would be injective if all~$A_i$, $i\in I$, were
mutually disjoint.) We shall call $\mu,\,\mu'\in\mathfrak M(\mathcal
A)$ {\it equivalent in\/}~$\mathfrak M(\mathcal A)$, and write
$\mu\cong\mu'$, whenever their $R$-images coincide.

It follows from~(\ref{R}) that, for given $\mu,\,\mu_1\in\mathfrak
M(\mathcal A)$ and $x\in\mathbf X$,
\begin{align*}
\kappa(x,R\mu)&=\sum_{i\in I}\,\alpha_i\kappa(x,\mu^i),\\[4pt]
\kappa(R\mu,R\mu_1)&=\sum_{i,j\in
I}\,\alpha_i\alpha_j\kappa(\mu^i,\mu_1^j),
\end{align*}
each of the above identities being understood in the sense that any
of its sides is well defined whenever so is the other, and then they
coincide. We shall call
\begin{align*}
\kappa(x,\mu)&:=\kappa(x,R\mu)\\
\intertext{the value of the {\it potential\/} of~$\mu$ at a point
$x$, and} \kappa(\mu,\mu_1)&:=\kappa(R\mu,R\mu_1)
\end{align*}
the {\it mutual energy\/} of~$\mu$ and~$\mu_1$~--- of course,
provided the right-hand side of the corresponding relation is well
defined. For $\mu\equiv\mu_1$ we get the {\it energy\/}
$\kappa(\mu,\mu)$ of~$\mu$.

Since we make no difference between $\mu\in\mathfrak M(\mathcal A)$
and $R\mu$ when dealing with their energies or potentials, we shall
sometimes call a measure associated with a condenser simply a {\it
measure\/}~--- certainly, if this causes no confusion.

Let $\mathcal E(\mathcal A)$ consist of all $\mu\in\mathfrak
M(\mathcal A)$ with finite energy $\kappa(\mu,\mu)=:\|\mu\|^2$. Then
$\mathcal E(\mathcal A)$ is {\it convex\/} and can be treated as a
semimetric space with the semimetric
\begin{equation}
\|\mu_1-\mu_2\|:=\|R\mu_1-R\mu_2\|,\quad\mu_1,\,\mu_2\in\mathcal
E(\mathcal A);\label{seminorm}
\end{equation}
the topology on $\mathcal E(\mathcal A)$ defined by means of this
semimetric will be called {\it strong\/}.

Two elements of $\mathcal E(\mathcal A)$, $\mu_1$ and~$\mu_2$, are
called {\it equivalent in\/}~$\mathcal E(\mathcal A)$ if
$\|\mu_1-\mu_2\|=0$. If, in addition, $\kappa$ is assumed to be
strictly positive definite, then the equivalence in~$\mathcal
E(\mathcal A)$ implies that in~$\mathfrak M(\mathcal A)$, namely
then $\mu_1\cong\mu_2$.\medskip

{\bf 3.4.} For measures associated with a condenser, it is also
reasonable to introduce the following concept of convergence,
actually corresponding to the vague convergence by coordinates. Let
$S$ denote a directed set of indices, and let $\mu_s$, $s\in S$, and
$\mu_0$ be given elements of the class~$\mathfrak
M(\,\overline{\mathcal A}\,)$.\medskip

{\bf Definition 3.} A net $(\mu_s)_{s\in S}$ is said to converge to
$\mu_0$ $\mathcal A$-{\it va\-gue\-ly\/} if
$$
\mu^i_s\to\mu_0^i\quad\mbox{vaguely for all \ } i\in I.
$$

Since $\mathfrak M(\mathbf X)$ is a Hausdorff space, an $\mathcal
A$-vague limit in~$\mathfrak M(\,\overline{\mathcal A}\,)$ is unique
(if exists).\medskip

{\bf Remark 3.} The $\mathcal A$-vague convergence of $(\mu_s)_{s\in
S}$ to $\mu_0$ certainly implies the vague convergence of
$(R\mu_s)_{s\in S}$ to~$R\mu_0$. By using the Tietze-Urysohn
extension theorem (see, e.\,g.,~\cite[Th.~0.2.13]{E2}), one can see
that, for the converse to be true, it is necessary and sufficient
that all $\overline{A_i}$, $i\in I$, be mutually disjoint.\bigskip

{\bf\large 4. Interior capacities of
$(m,p)$-condensers}\nopagebreak\medskip

{\bf 4.1.} Let $\mathcal H$ be a set in the pre-Hilbert
space~$\mathcal E$ or in the semimetric space~$\mathcal E(\mathcal
A)$, where $\mathcal A$ is a given $(m,p)$-con\-den\-ser. In either
case, let us introduce the quantity
$$
\|\mathcal H\|^2:=\inf_{\nu\in\mathcal H}\,\|\nu\|^2,
$$
interpreted as $+\infty$ if $\mathcal H$ is empty. If $\|\mathcal
H\|^2<\infty$, one can consider the variational problem on the
existence of $\lambda=\lambda(\mathcal H)\in\mathcal H$ with minimal
energy
$$
\|\lambda\|^2=\|\mathcal H\|^2;
$$
such a problem will be referred to as the $\mathcal H$-{\it
problem\/}. The $\mathcal H$-problem is said to be {\it solvable\/}
if a minimizer~$\lambda(\mathcal H)$ exists.

The following elementary lemma is a slight generalization of
\cite[Lemma~4.1.1]{F1}.\medskip

{\bf Lemma 1.} {\it Suppose $\mathcal H$ is convex, and
$\lambda=\lambda(\mathcal H)$ exists. Then for any $\nu\in\mathcal
H$,}
\begin{equation}
\|\nu-\lambda\|^2\leqslant\|\nu\|^2-\|\lambda\|^2.\label{lemma1}
\end{equation}

{\bf Proof.} Assume $\mathcal H\subset\mathcal E$. For every
$t\in[0,\,1]$, the measure $\mu:=(1-t)\lambda+t\nu$ belongs
to~$\mathcal H$, and therefore $\|\mu\|^2\geqslant\|\lambda\|^2$.
Evaluating~$\|\mu\|^2$ and then letting $t$ tend to zero, we get
$\kappa(\nu,\lambda)\geqslant\|\lambda\|^2$, and (\ref{lemma1})
follows (see~\cite{F1}).

Suppose now $\mathcal H\subset\mathcal E(\mathcal A)$. Then
$R\mathcal H:=\{R\nu:\ \nu\in\mathcal H\}$ is a convex subset
of~$\mathcal E$, while $R\lambda$ is a minimizer in the $R\mathcal
H$-problem. What has been shown thus yields
$$\|R\nu-R\lambda\|^2\leqslant\|R\nu\|^2-\|R\lambda\|^2,$$
which gives (\ref{lemma1}) when combined
with~(\ref{seminorm}).\medskip

We shall be concerned with the $\mathcal H$-problem for various
specific~$\mathcal H$ related to the notion of {\it interior
capacity\/} of an $(m,p)$-con\-den\-ser (in particular, of a set);
see~Sec.~4.2 and Sec.~6 below for their definitions.\medskip

{\bf 4.2.} Fix a continuous function $g:\mathbf X\to(0,\infty)$ and
a numerical vector $a=(a_i)_{i\in I}$ with $a_i>0$, $i\in I$. Given
a kernel~$\kappa$ and an $(m,p)$-condenser~$\mathcal A$ in~$\mathbf
X$, write
$$
\mathcal E(\mathcal A,a,g):=\Bigl\{\mu\in\mathcal E(\mathcal A):  \ \int
g\,d\mu^i=a_i\mbox{ \ for all \ }i\in I\Bigr\}.
$$

{\bf Definition 4.} We shall call the value
\begin{equation}
{\rm cap}\,\mathcal A:={\rm cap}\,(\mathcal
A,a,g):=\frac{1}{\|\mathcal E(\mathcal A,a,g)\|^2} \label{def}
\end{equation}
the ({\it interior}) {\it capacity\/} of an $(m,p)$-condenser
$\mathcal A$ (with respect to~$\kappa$, $a$, and~$g$).\medskip

Here and in the sequel, we adopt the convention that $1/0=+\infty$.
It follows immediately from the positive definiteness of the kernel
that
$$
0\leqslant{\rm cap}\,(\mathcal A,a,g)\leqslant\infty.
$$

{\bf Remark 4.} If $I$ is a singleton, any $(m,p)$-con\-den\-ser
consists of just one set, $A_1$, positively signed. If moreover
$g=1$ and $a_1=1$, then the notion of interior capacity of a
condenser, defined above, certainly reduces to the notion of
interior capacity of a set (see~\cite{F1}). We denote it
by~$C(\,\cdot\,)$ as well.\medskip

{\bf Remark 5.} In the case of the Newton kernel in $\mathbb R^3$,
the notion of capacity of a condenser~$\mathcal A$ has an evident
electrostatic interpretation. In the framework of the corresponding
electrostatics problem, the function~$g$ serves as a characteristic
of nonhomogeneity of the conductors $A_i$, $i\in I$.\medskip

{\bf 4.3.} On $\mathfrak C_{m,p}=\mathfrak C_{m,p}(\mathbf X)$, it
is natural to introduce the ordering relation~$\prec$ by declaring
$\mathcal A'\prec\mathcal A$ to mean that $A_i'\subset A_i$ for all
$i\in I$. Here, $\mathcal A'=(A_i')_{i\in I}$. Then ${\rm
cap}\,(\,\cdot\,,a,g)$ is a nondecreasing function of a condenser,
namely
\begin{equation}
{\rm cap}\,(\mathcal A',a,g)\leqslant{\rm cap}\,(\mathcal
A,a,g)\quad\mbox{whenever \ }\mathcal A'\prec\mathcal A.
\label{increas}
\end{equation}

Given $\mathcal A\in\mathfrak C_{m,p}$, denote by $\{\mathcal
K\}_\mathcal A$ the increasing ordered family of all compact
condensers $\mathcal K=(K_i)_{i\in I}\in\mathfrak C_{m,p}$ such that
$\mathcal K\prec\mathcal A$.\medskip

{\bf Lemma 2.} {\it If $\mathcal K$ ranges over $\{\mathcal
K\}_\mathcal A$, then}
\begin{equation}
{\rm cap}\,(\mathcal A,a,g)=\lim_{\mathcal K\uparrow\mathcal A}\,
{\rm cap}\,(\mathcal K,a,g).\label{contin'}
\end{equation}

{\bf Proof.} We can certainly assume ${\rm cap}\,(\mathcal A,a,g)$
to be nonzero, since otherwise the lemma follows at once
from~(\ref{increas}). Then the set $\mathcal E(\mathcal A,a,g)$ must
be nonempty; fix~$\mu$, one of its elements. For any $\mathcal
K\in\{\mathcal K\}_\mathcal A$ and $i\in I$, let $\mu_\mathcal K^i$
denote the trace of~$\mu^i$ upon~$K_i$. Applying Lemma~1.2.2
from~\cite{F1}, we get
\begin{align}
\int g\,d\mu^i&=\lim_{\mathcal K\uparrow\mathcal A}\,\int
g\,d\mu_\mathcal K^i,\quad i\in I,\label{mon}\\
\kappa(\mu^i,\mu^j)&=\lim_{\mathcal K\uparrow\mathcal
A}\,\kappa(\mu_\mathcal K^i,\mu_\mathcal K^j),\quad i,\,j\in I.
\label{mon'}
\end{align}
Thus for $\mathcal K\in\{\mathcal K\}_\mathcal A$ large enough,
$\int g\,d\mu_\mathcal K^i\ne0$ for all $i\in I$, and consequently
$$
\sum_{i\in I}\,\frac{\alpha_ia_i}{\int g\,d\mu_\mathcal
K^i}\,\mu_\mathcal K^i\in\mathcal E(\mathcal K,a,g).
$$
Together with (\ref{mon}) and (\ref{mon'}), this yields
$$
\|\mu\|^2=\lim_{\mathcal K\uparrow\mathcal A}\, \sum_{i,j\in I}\,
\kappa\Bigl(\frac{\alpha_ia_i}{\int g\,d\mu_\mathcal
K^i}\,\mu_\mathcal K^i, \frac{\alpha_ja_j}{\int g\,d\mu_\mathcal
K^j}\,\mu_\mathcal K^j\Bigr) \geqslant\lim_{\mathcal
K\uparrow\mathcal A}\,\|\mathcal E(\mathcal K,a,g)\|^2,
$$
and hence, in view of the arbitrary choice of $\mu\in\mathcal
E(\mathcal A,a,g)$,
$$
\|\mathcal E(\mathcal A,a,g)\|^2\geqslant\lim_{\mathcal
K\uparrow\mathcal A}\,\|\mathcal E(\mathcal K,a,g)\|^2.
$$
Since the converse inequality is obvious from~(\ref{increas}), the
proof is complete.\medskip

Let $\mathcal E^0(\mathcal A,a,g)$ denote the class of all
$\mu\in\mathcal E(\mathcal A,a,g)$ such that, for every $i\in I$,
the support~$S(\mu^i)$ of~$\mu^i$ is compact and contained
in~$A_i$.\medskip

{\bf Corollary 1.} {\it The capacity ${\rm cap}\,(\mathcal A,a,g)$
remains unchanged if the class $\mathcal E(\mathcal A,a,g)$ in its
definition is replaced by $\mathcal E^0(\mathcal A,a,g)$.  In other
words,}
$$
\|\mathcal E(\mathcal A,a,g)\|^2=\|\mathcal E^0(\mathcal A,a,g)\|^2.
$$

{\bf Proof.} We can certainly assume ${\rm cap}\,\mathcal A$ to be
nonzero, since otherwise the corollary follows immediately from the
inclusion $\mathcal E^0(\mathcal A,a,g)\subset\mathcal E(\mathcal
A,a,g)$. Then, by~(\ref{increas}) and~(\ref{contin'}), for every
$\varepsilon>0$ there exists a compact condenser $\mathcal
K\prec\mathcal A$ such that
$$\|\mathcal E(\mathcal K,a,g)\|^2\leqslant\|\mathcal E(\mathcal A,a,g)\|^2+\varepsilon.$$
This leads to the claimed assertion when combined with the relation
$$\|\mathcal E(\mathcal K,a,g)\|^2\geqslant\|\mathcal
E^0(\mathcal A,a,g)\|^2\geqslant\|\mathcal E(\mathcal A,a,g)\|^2.$$

{\bf 4.4.} Unless explicitly stated otherwise, in all that follows
it is assumed that
\begin{equation}
{\rm cap}\,(\mathcal A,a,g)>0.\label{nonzero1}
\end{equation}

{\bf Lemma 3.} {\it The assumption\/} (\ref{nonzero1}) {\it is
equivalent to the following one\/}:
\begin{equation}
C(A_i)>0\quad\mbox{\it for all \ } i\in I. \label{cap1}
\end{equation}

{\bf Proof.} Indeed, ${\rm cap}\,(\mathcal A,a,g)$ is nonzero if and
only if $\mathcal E(\mathcal A,a,g)$ is nonempty. As $g$ is
positive, for the latter to happen, it is necessary and sufficient
that, for every $i\in I$, there exists a nonzero nonnegative measure
of finite energy concentrated on~$A_i$. Since this is equivalent
to~(\ref{cap1}) by~\cite[Lemma~2.3.1]{F1}, the proof is
complete.\medskip

In the following assertion, providing necessary and sufficient
conditions for ${\rm cap}\,\mathcal A$ to be finite, we assume
$g|_A$ to have a strictly positive lower bound (say~$L$).\medskip

{\bf Lemma 4.} {\it For ${\rm cap}\,(\mathcal A,a,g)$ to be finite,
it is necessary that
\begin{equation}
C(A_j)<\infty\quad\mbox{for some \ } j\in I.\label{j}
\end{equation}
This condition is also sufficient if it is additionally assumed that
$\mathcal A$ is closed, $g|_A$ bounded, and $\kappa$ bounded from
above on~$A^+\times A^-$ and perfect.}\medskip

{\bf Proof.} Let ${\rm cap}\,\mathcal A<\infty$, and assume, on the
contrary, that
\begin{equation}
C(A_i)=\infty\quad\mbox{for all \ } i\in I.\label{lemma6}
\end{equation}
Then, for every~$i$, there exist probability measures
$\nu_n^i\in\mathcal E^+(A_i)$, $n\in\mathbb N$, of compact support
such that
$$
\|\nu_n^i\|\to0\quad(n\to\infty).
$$
Since
$$
\mu_n:=\sum_{i\in I}\,\frac{\alpha_i a_i\nu_n^i}{\int
g\,d\nu_n^i}\in\mathcal E(\mathcal A,a,g),\quad n\in\mathbb N,
$$
and
$$
\|\mu_n\|\leqslant L^{-1}\,\sum_{i\in I}\,a_i\|\nu_n^i\|,
$$
we arrive at a contradiction by letting $n$ tend to $\infty$.

Assume now all the conditions of the remaining part of the lemma to
be satisfied, and let (\ref{j}) be true. Then,
by~\cite[Lemma~13]{Z5}, there exists $\zeta\in\mathcal E(\mathcal
A)$ with the properties that $\int g\,d\zeta^j=a_j$ (hence,
$\zeta\not\equiv0$) and
$$\|\zeta\|^2=\|\mathcal E(\mathcal A,a,g)\|^2.$$
Since $\kappa$ is strictly positive definite, this yields ${\rm
cap}\,\mathcal A<\infty$, as was to be proved.\medskip

{\bf 4.5.} Because of~(\ref{nonzero1}), we are naturally led to the
$\mathcal E(\mathcal A,a,g)$-{\it problem\/} (cf.~Sec.~4.1), i.\,e.,
the problem on the existence of $\lambda\in\mathcal E(\mathcal
A,a,g)$ with minimal energy
$$
\|\lambda\|^2=\|\mathcal E(\mathcal A,a,g)\|^2;
$$
the $\mathcal E(\mathcal A,a,g)$-problem might certainly be regarded
as the main minimum-problem of the theory of interior capacities of
condensers. The collection (possibly empty) of all minimizing
measures~$\lambda$ in this problem will be denoted by $\mathcal
S(\mathcal A,a,g)$.

If moreover ${\rm cap}\,(\mathcal A,a,g)$ is finite, let us look, as
well, at the $\mathcal E(\mathcal A,a\,{\rm cap}\,\mathcal
A,g)$-{\it prob\-lem\/}. By reasons of homogeneity, both the
$\mathcal E(\mathcal A,a,g)$- and the $\mathcal E(\mathcal A,a\,{\rm
cap}\,\mathcal A,g)$-problems are simultaneously either solvable or
unsolvable, and their extremal values are related to each other by
the following law:
\begin{equation}
\frac{1}{\|\mathcal E(\mathcal A,a,g)\|^2}=\|\mathcal E(\mathcal
A,a\,{\rm cap}\,\mathcal A,g)\|^2.
\label{iden}
\end{equation}

Assume for a moment that $\mathcal A$ is compact. Since the mapping
$$\nu\mapsto\int g\,d\nu,\quad\nu\in\mathfrak M^+(K),$$
where $K\subset\mathbf X$ is a compact set, is vaguely continuous,
$\mathcal E(\mathcal A,a,g)$ is compact in the $\mathcal A$-vague
topology. Therefore, if $\kappa$ is additionally assumed to be
continuous on~$A^+\times A^-$ (which, due to~(\ref{non}), is always
the case for either of the classical kernels), then the
energy~$\|\mu\|^2$ is $\mathcal A$-vaguely lower semicontinuous
on~$\mathcal E(\mathcal A)$, and the solvability of both the
problems immediately follows (cf.~\cite[Th.~2.6]{O}).

But if $\mathcal A$ is {\it noncompact\/}, then the class $\mathcal
E(\mathcal A,a,g)$ is no longer $\mathcal A$-vaguely compact and the
problems become quite nontrivial. Moreover, it has recently been
shown by the author that, in the noncompact case, the problems are
in general {\it unsolvable\/} and this phenomenon occurs even under
very natural assumptions (e.\,g., for the Newton, Green, or Riesz
kernels in~$\mathbb R^n$, $n\geqslant 2$, and closed condensers).

In particular, it was proved in~\cite{Z4} that, if $\mathcal A$ is
closed, $\kappa$~is perfect, and bounded and continuous
on~$A^+\times A^-$, and satisfies the generalized maximum principle
(see, e.\,g.,~\cite[Chap.~VI]{L}), while $g|_A$ is bounded and has a
strictly positive lower bound, then for either of the $\mathcal
E(\mathcal A,a,g)$- and the $\mathcal E(\mathcal A,a\,{\rm
cap}\,\mathcal A,g)$-prob\-lems to be solvable for any vector~$a$,
it is necessary and sufficient that
$$
C(A_i)<\infty\quad\mbox{for all \ } i\in I.
$$
If moreover there exists $i_0\in I$ such that
$$
C(A_{i_0})=\infty,
$$
then both the problems are unsolvable for all $a=(a_i)_{i\in I}$
with $a_{i_0}$ large enough.

In~\cite[Th.~1]{Z6}, the last statement was sharpened. It was shown
that if, in addition to all the preceding assumptions, for all
$i\neq i_0$,
$$
C(A_i)<\infty\quad\mbox{and}\quad A_i\cap A_{i_0}=\varnothing,
$$
while $\kappa(\cdot,y)\to0$ (as $y\to\infty$) uniformly on compact
sets, then there exists a number $\Lambda_{i_0}\in[0,\infty)$ such
that the problems are unsolvable if and only if
$$a_{i_0}>\Lambda_{i_0}.$$

{\bf Remark 6.} It was actually shown in \cite{Z6} that
$$
\Lambda_{i_0}=\int g\,d\tilde{\lambda}^{i_0},
$$
where $\tilde{\lambda}$ is a minimizer (it exists) in the auxiliary
$\mathcal H$-problem for
$$
\mathcal H:=\Bigl\{\mu\in\mathcal E(\mathcal A):\quad\int
g\,d\mu^i=a_i\mbox{ \ for all \ } i\neq i_0\Bigr\}.
$$

{\bf Remark 7.} The mentioned results were actually obtained
in~\cite{Z4,Z6} for the energy evaluated in the presence of an
external field.\medskip

{\bf 4.6.} In view of the results reviewed in Sec.~4.5, it was
particularly interesting to find statements of variational problems
{\it dual\/} to the $\mathcal E(\mathcal A,a\,{\rm cap}\,\mathcal
A,g)$-problem (and hence providing some new {\it equivalent\/}
definitions of ${\rm cap}\,\mathcal A$), but {\it solvable\/} for
any condenser~$\mathcal A$ (e.\,g., even nonclosed) and any
vector~$a$. We have succeeded in this under the following
conditions, which will always be tacitly assumed.

From now on, in addition to~(\ref{nonzero1}), the following {\bf
standing assumptions} are always required: $\kappa$~is consistent,
and either
$$
I^-=\varnothing\quad\mbox{(i.\,e., \ $p=0$)},
$$
or both the conditions are satisfied
\begin{align}
g_{\min}:=\inf_{x\in A}\,g(x)&>0,\label{g}\\[5pt]
\sup_{x\in A^+,\ y\in A^-}\,\kappa(x,y)&<\infty.\label{bou}
\end{align}

{\bf Remark 8.} These assumptions on a kernel are not too
restrictive. In particular, they all are satisfied by the Newton,
Riesz, or Green kernels in~$\mathbb R^n$ provided the Euclidean
distance between the opposite-signed plates of a condenser is
nonzero.\bigskip

{\bf\large 5. $\mathcal A$-vague and strong cluster sets of
minimizing nets}\nopagebreak\medskip

To formulate the results obtained, we shall need the following
notation.

{\bf 5.1.} Denote by $\mathbb M(\mathcal A,a,g)$ the class of all
$(\mu_t)_{t\in T}\subset\mathcal E^0(\mathcal A,a,g)$ such that
\begin{equation}
\lim_{t\in T}\,\|\mu_t\|^2=\|\mathcal E(\mathcal A,a,g)\|^2.
\label{min}
\end{equation}
This class is not empty, which is clear from~(\ref{nonzero1}) in
view of Corollary~1.

Let $\mathcal M(\mathcal A,a,g)$ (respectively, $\mathcal
M'(\mathcal A,a,g)$) consist of all limit points of the nets
$(\mu_t)_{t\in T}\in\mathbb M(\mathcal A,a,g)$ in the $\mathcal
A$-vague topology of the space~$\mathfrak M(\,\overline{\mathcal
A}\,)$ (respectively, in the strong topology of the semimetric space
$\mathcal E(\,\overline{\mathcal A}\,)$). Also write
$$
\mathcal E(\mathcal A,\leqslant\!a,g):=\Bigl\{\mu\in\mathcal
E(\mathcal A):\quad\int g\,d\mu^i\leqslant a_i\mbox{ \ for all \ }
i\in I\Bigr\}.
$$

With the preceding notation and under our standing assumptions
(see~Sec.~4.6), there holds the following lemma, to be proved in
Sec.~11 below.\medskip

{\bf Lemma 5.} {\it Given $(\mu_t)_{t\in T}\in\mathbb M(\mathcal
A,a,g)$, there exist its $\mathcal A$-vague cluster points; hence,
the class $\mathcal M(\mathcal A,a,g)$ is nonempty. Moreover,
\begin{equation}
\mathcal M(\mathcal A,a,g)\subset\mathcal M'(\mathcal
A,a,g)\cap\mathcal E(\,\overline{\mathcal
A},\leqslant\!a,g).\label{MMprime}
\end{equation}
Furthermore, for every $\chi\in\mathcal M'(\mathcal A,a,g)$,
\begin{equation}
\lim_{t\in T}\,\|\mu_t-\chi\|^2=0,\label{strongly}
\end{equation}
and hence $\mathcal M'(\mathcal A,a,g)$ forms an equivalence class
in~$\mathcal E(\,\overline{\mathcal A}\,)$.}\medskip

It follows from (\ref{min})\,--\,(\ref{strongly}) that
$$
\|\zeta\|^2=\|\mathcal E(\mathcal A,a,g)\|^2\quad\mbox{for all \
}\zeta\in\mathcal M(\mathcal A,a,g).
$$
Also observe that, if $\mathcal A=\mathcal K$ is compact, then
moreover $\mathcal M(\mathcal K,a,g)\subset\mathcal E(\mathcal
K,a,g)$, which together with the preceding relation proves the
following assertion.\medskip

{\bf Corollary 2.} {\it If $\mathcal A=\mathcal K$ is compact, then
the $\mathcal E(\mathcal K,a,g)$-problem is solvable. Actually,}
\begin{equation}
\mathcal S(\mathcal K,a,g)=\mathcal M(\mathcal K,a,g).\label{S}
\end{equation}

{\bf 5.2.} When approaching $\mathcal A$ by the increasing
family~$\{\mathcal K\}_\mathcal A$ of the compact condensers
$\mathcal K\prec\mathcal A$, we shall always suppose all
those~$\mathcal K$ to be of capacity nonzero. This involves no loss
of generality, which is clear from~(\ref{nonzero1}) and Lemma~2.

Then Corollary 2 enables us to introduce the (nonempty) class
$\mathbb M_0(\mathcal A,a,g)$ of all nets $(\lambda_\mathcal
K)_{\mathcal K\in\{\mathcal K\}_\mathcal A}$, where
$\lambda_\mathcal K\in\mathcal S(\mathcal K,a,g)$ is arbitrarily
chosen. Let $\mathcal M_0(\mathcal A,a,g)$ consist of all $\mathcal
A$-vague cluster points of those nets. Since, by Lemma~2,
$$\mathbb M_0(\mathcal A,a,g)\subset\mathbb M(\mathcal A,a,g),$$
application of Lemma~5 yields the following assertion.\medskip

{\bf Corollary 3.} {\it The class $\mathcal M_0(\mathcal A,a,g)$ is
nonempty, and}
$$
\mathcal M_0(\mathcal A,a,g)\subset\mathcal M(\mathcal
A,a,g)\subset\mathcal M'(\mathcal A,a,g).
$$

{\bf Remark 9.} Each of the cluster sets, $\mathcal M_0(\mathcal
A,a,g)$, $\mathcal M(\mathcal A,a,g)$ and $\mathcal M'(\mathcal
A,a,g)$, plays an important role in our study. However, if $\kappa$
is additionally assumed to be strictly positive definite (hence,
perfect), while $\overline{A_i}$, $i\in I$, are mutually disjoint,
then all these classes coincide and consist of just one
element.\medskip

{\bf 5.3.} Also the following notation will be required. Given
$\chi\in\mathcal M'(\mathcal A,a,g)$, write
$$
\mathcal M'_\mathcal E(\mathcal A,a,g):=\bigl[R\chi\bigr]_\mathcal
E\,.
$$
This equivalence class does not depend on the choice of $\chi$,
which is clear from Lemma~5. Lemma~5 also yields that, for any
$(\mu_t)_{t\in T}\in\mathbb M(\mathcal A,a,g)$ and any
$\nu\in\mathcal M'_\mathcal E(\mathcal A,a,g)$, $R\mu_t\to\nu$ in
the strong topology of the pre-Hilbert space~$\mathcal E$.\bigskip

{\bf\large 6. Extremal problems dual to the main minimum-pro\-blem
of the theory of interior capacities of
condensers}\nopagebreak\medskip

Throughout Sec.~6, as usual, we are keeping all our standing
assumptions, stated in~Sec.~4.6.\medskip

{\bf 6.1.} A proposition $R(x)$ involving a variable point
$x\in\mathbf X$ is said to subsist {\it nearly
everywhere\/}~(n.\,e.) in~$E$, where $E$ is a given subset
of~$\mathbf X$, if the set of all $x\in E$ for which $R(x)$ fails to
hold is of interior capacity zero. See,~e.\,g.,~\cite{F1}.

If $C(E)>0$ and $f$ is a universally measurable function bounded
from below nearly everywhere in~$E$, write
$$
"\!\inf_{x\in E}\!"\,f(x):=\sup\,\bigl\{q: \ f(x)\geqslant
q\quad\mbox{n.\,e.~in \ } E\bigr\}.
$$
Then
$$
f(x)\geqslant"\!\inf_{x\in E}\!"\,f(x)\quad\mbox{n.\,e.~in \ } E,
$$
which is seen from the fact that the interior
capacity~$C(\,\cdot\,)$ is countably subadditive on sets $U_n\cap
E$, $n\in\mathbb N$, where $U_n$ are universally measurable, whereas
$E$ is arbitrary (see Lemma~2.3.5 in~\cite{F1} and the remark
attached to it).\medskip

{\bf 6.2.} Let $\hat{\Gamma}=\hat{\Gamma}(\mathcal A,a,g)$ denote
the class of all Radon measures $\nu\in\mathcal E$ such that there
exist real numbers $c_i(\nu)$, $i\in I$, satisfying the relations
\begin{equation}
\alpha_i a_i\kappa(x,\nu)\geqslant c_i(\nu)g(x)\quad\mbox{n.\,e. in
\ } A_i,\quad i\in I, \label{adm1}
\end{equation}
\begin{equation}
\sum_{i\in I}\,c_i(\nu)\geqslant1. \label{adm2}
\end{equation}
The property of subadditivity of~$C(\,\cdot\,)$, mentioned above,
implies that $\hat{\Gamma}$ is {\it convex\/}.

The following assertion, to be proved in Sec.~14 below, holds
true.\medskip

{\bf Theorem 2.} {\it Under the standing assumptions,}
\begin{equation}
\|\hat{\Gamma}(\mathcal A,a,g)\|^2={\rm cap}\,(\mathcal A,a,g).
\label{id}
\end{equation}

If $\|\hat{\Gamma}(\mathcal A,a,g)\|^2<\infty$, we are interested in
the $\hat{\Gamma}(\mathcal A,a,g)$-{\it problem\/} (cf.~Sec.~4.1),
i.\,e., the problem on the existence of
$\hat{\omega}\in\hat{\Gamma}(\mathcal A,a,g)$ with minimal energy
$$
\|\hat{\omega}\|^2=\|\hat{\Gamma}(\mathcal A,a,g)\|^2;
$$
the collection of all those $\hat\omega$ will be denoted by
$\hat{\mathcal G}=\hat{\mathcal G}(\mathcal A,a,g)$.

A minimizing measure $\hat{\omega}$ can be shown to be {\it
unique\/} up to a summand of seminorm zero (and, hence, it is unique
whenever the kernel under consideration is strictly positive
definite). Actually, the following stronger result holds
true.\medskip

{\bf Lemma 6.} {\it If $\hat\omega$ exists, $\hat{\mathcal
G}(\mathcal A,a,g)$ forms an equivalence class in~$\mathcal
E$.}\medskip

{\bf Proof.} Since $\hat{\Gamma}$ is convex, Lemma~1 yields that
$\hat{\mathcal G}$ is contained in an equivalence class in~$\mathcal
E$. To prove that $\hat{\mathcal G}$ actually coincides with that
equivalence class, it suffices to show that, if $\nu$ belongs
to~$\hat{\Gamma}$, then so do all measures equivalent to~$\nu$
in~$\mathcal E$. But this follows at once from the property of
subadditivity of~$C(\,\cdot\,)$, mentioned above, and the fact that
the potentials of any two equivalent in~$\mathcal E$ measures
coincide nearly everywhere in~$\mathbf
X$~\cite[Lemma~3.2.1]{F1}.\medskip

{\bf 6.3.} Assume for a moment that ${\rm cap}\,(\mathcal A,a,g)$ is
finite. When combined with~(\ref{def}) and~(\ref{iden}), Theorem~2
shows that the $\hat{\Gamma}(\mathcal A,a,g)$-pro\-blem and, on the
other hand, the $\mathcal E(\mathcal A,a\,{\rm cap}\,\mathcal
A,g)$-pro\-blem have the same infimum, equal to the capacity ${\rm
cap}\,\mathcal A$, and so these two variational problems are {\it
dual}.

But what is surprising is that their infimum, ${\rm cap}\,\mathcal
A$, turns out to be always an actual minimum in the former extremal
problem, while this is not the case for the latter one
(see~Sec.~4.5). In fact, the following statement on the solvability
of the $\hat{\Gamma}(\mathcal A,a,g)$-problem, to be proved in
Sec.~14 below, holds true.\medskip

{\bf Theorem 3.} {\it Under the standing assumptions, if moreover
${\rm cap}\,\mathcal A<\infty$, then the class $\hat{\mathcal
G}(\mathcal A,a,g)$ is nonempty and can be given by the formula
\begin{equation}
\hat{\mathcal G}(\mathcal A,a,g)=\mathcal M'_\mathcal E(\mathcal
A,a\,{\rm cap}\,\mathcal A,g).\label{desc}
\end{equation}
The numbers $c_i(\hat{\omega})$, $i\in I$, satisfying
both\/}~(\ref{adm1}) {\it and\/}~(\ref{adm2}) {\it for
$\hat{\omega}\in\hat{\mathcal G}(\mathcal A,a,g)$, are determined
uniquely, do not depend on the choice of~$\hat{\omega}$, and can be
written in either of the forms
\begin{align}
c_i(\hat{\omega})&=\alpha_i\,{\rm cap}\,\mathcal
A^{-1}\kappa(\zeta^i,\zeta),\label{snt1}\\[4pt]
c_i(\hat{\omega})&=\alpha_i\,{\rm cap}\,\mathcal A^{-1}\lim_{s\in
S}\,\kappa(\mu_s^i,\mu_s),\label{cnt}
\end{align}
where $\zeta\in\mathcal M(\mathcal A,a\,{\rm cap}\,\mathcal A,g)$
and $(\mu_s)_{s\in S}\in\mathbb M(\mathcal A,a\,{\rm cap}\,\mathcal
A,g)$ are arbitrarily given.}\medskip

The following two assertions, providing additional information
about~$c_i(\hat{\omega})$, $i\in I$, can be obtained directly from
the preceding theorem.\medskip

{\bf Corollary 4.} {\it Given $\hat{\omega}\in\hat{\mathcal
G}(\mathcal A,a,g)$, it follows that}
\begin{equation}
c_i(\hat{\omega})="\!\inf_{x\in A_i}\!"\,\,\frac{\alpha_i
a_i\kappa(x,\hat{\omega})}{g(x)},\quad i\in
I\smallskip.\label{ess''}
\end{equation}

{\bf Corollary 5.} {\it The inequality\/}~(\ref{adm2}) {\it for
$\hat{\omega}\in\hat{\mathcal G}(\mathcal A,a,g)$ is actually an
equality; i.\,e.} \begin{equation}\sum_{i\in
I}\,c_i(\hat\omega)=1.\label{1''}\end{equation}

{\bf Remark 10.} Assume for a moment that ${\rm cap}\,\mathcal A=0$.
Then, by Lemma~3, there exists $i\in I$ (say $i=1$) with $C(A_i)=0$.
Hence, the measure $\nu_0=0$ belongs to $\hat{\Gamma}(\mathcal
A,a,g)$ since it satisfies both~(\ref{adm1}) and~(\ref{adm2}) with
$c_i(\nu_0)$, $i\in I$, where
$$c_1(\nu_0)\geqslant1\quad\mbox{and}\quad c_i(\nu_0)=0,\quad i\ne1.$$ This implies that
the identity~(\ref{id}) actually holds true in the degenerate case
${\rm cap}\,\mathcal A=0$ as well, and then $\hat{\mathcal
G}(\mathcal A,a,g)$ consists of all $\nu\in\mathcal E$ of seminorm
zero. What then, however, fails to hold is the statement on the
uniqueness of~$c_i(\hat{\omega})$.\medskip

{\bf 6.4.} Let $\hat{\Gamma}_*(\mathcal A,a,g)$ consist of all
$\nu\in\hat{\Gamma}(\mathcal A,a,g)$ for which the
inequality~(\ref{adm2}) is actually an equality. By arguments
similar to those that have been applied above, one can see that
$\hat{\Gamma}_*(\mathcal A,a,g)$ is convex, and hence all the
solutions to the minimal energy problem over this class form an
equivalence class in~$\mathcal E$. Combining this with Theorems~2,~3
and Corollary~5 leads to the following assertion.\medskip

{\bf Corollary 6.} {\it Under the standing assumptions,
$$
\|\hat{\Gamma}_*(\mathcal A,a,g)\|^2={\rm cap}\,(\mathcal A,a,g).
$$
If moreover ${\rm cap}\,\mathcal A<\infty$, then the
$\hat{\Gamma}_*(\mathcal A,a,g)$-problem is solvable and the class
$\hat{\mathcal G}_*(\mathcal A,a,g)$ of all its solutions is given
by the formula}
$$\hat{\mathcal G}_*(\mathcal A,a,g)=
\mathcal M'_\mathcal E(\mathcal A,a\,{\rm cap}\,\mathcal A,g).
$$

{\bf Remark 11.} Theorem~2 and Corollary~6 (cf.~also Theorem~4 and
Corollary~8 below) provide new equivalent definitions of the
capacity ${\rm cap}\,(\mathcal A,a,g)$. Note that, in contrast to
the initial definition (cf.~Sec.~4.2), no restrictions on the
supports and total masses of measures from the classes
$\hat{\Gamma}(\mathcal A,a,g)$ or $\hat{\Gamma}_*(\mathcal A,a,g)$
have been imposed; the only restriction involves their potentials.
These definitions of the capacity are actually new even in the
compact case; compare with~\cite{O}. They are not only of obvious
academic interest, but seem also to be important for numerical
computations.\medskip

{\bf 6.5.} Our next purpose is to formulate an $\mathcal H$-problem
such that it is still dual to the $\mathcal E(\mathcal A,a\,{\rm
cap}\,\mathcal A,g)$-problem and solvable, but now with $\mathcal H$
consisting of measures associated with a condenser.

Let $\Gamma(\mathcal A,a,g)$ consist of all $\mu\in\mathcal
E(\,\overline{\mathcal A}\,)$ for which both the
relations~(\ref{adm1}) and~(\ref{adm2}) hold (with~$\mu$ in place
of~$\nu$). In other words,
\begin{equation}
\Gamma(\mathcal A,a,g):=\left\{\mu\in\mathcal E(\,\overline{\mathcal
A}\,):\quad R\mu\in\hat{\Gamma}(\mathcal
A,a,g)\right\}.\label{Gamma}
\end{equation}
Observe that the class $\Gamma(\mathcal A,a,g)$ is {\it convex\/}
and
\begin{equation}
\|\Gamma(\mathcal A,a,g)\|^2\geqslant\|\hat{\Gamma}(\mathcal
A,a,g)\|^2. \label{GammahatGamma}
\end{equation}

We proceed to show that the inequality (\ref{GammahatGamma}) is
actually an equality, and that the minimal energy problem, if
considered over the class $\Gamma(\mathcal A,a,g)$, is still
solvable.\medskip

{\bf Theorem 4.} {\it Under the standing assumptions,
\begin{equation}
\|\Gamma(\mathcal A,a,g)\|^2={\rm cap}\,(\mathcal A,a,g).
\label{id2}
\end{equation}
If moreover ${\rm cap}\,(\mathcal A,a,g)<\infty$, then the
$\Gamma(\mathcal A,a,g)$-problem is solvable and the class $\mathcal
G(\mathcal A,a,g)$ of all its solutions~$\omega$ is given by the
formula}
\begin{equation}
\mathcal G(\mathcal A,a,g)=\mathcal M'(\mathcal A,a\,{\rm
cap}\,\mathcal A,g).\label{desc2}
\end{equation}

{\bf Proof.} We can certainly assume ${\rm cap}\,\mathcal A$ to be
finite, for if not, (\ref{id2})~is obtained directly from~(\ref{id})
and~(\ref{GammahatGamma}). Then, according to Lemma~5 with~$a\,{\rm
cap}\,\mathcal A$ instead of~$a$, the class $\mathcal M'(\mathcal
A,a\,{\rm cap}\,\mathcal A,g)$ is nonempty; fix~$\chi$, one of its
elements. It is clear from its definition and the
identity~(\ref{desc}) that $\chi\in\mathcal E(\,\overline{\mathcal
A}\,)$ and $R\chi\in\hat{\mathcal G}(\mathcal A,a,g)$. Hence,
by~(\ref{Gamma}), $\chi\in\Gamma(\mathcal A,a,g)$, and therefore
$$\|\hat{\Gamma}(\mathcal A,a,g)\|^2=\|\chi\|^2\geqslant\|\Gamma(\mathcal A,a,g)\|^2.$$
In view of~(\ref{id}) and~(\ref{GammahatGamma}), this
proves~(\ref{id2}) and, as well, the inclusion
$$\mathcal M'(\mathcal A,a\,{\rm cap}\,\mathcal A,g)\subset\mathcal G(\mathcal A,a,g).$$
But the right-hand side of this inclusion is an equivalence class
in~$\mathcal E(\,\overline{\mathcal A}\,)$, which is proved by the
convexity of $\Gamma(\mathcal A,a,g)$ and Lemma~1 in the same manner
as in the proof of Lemma~6. Since, by Lemma~5, also the left-hand
side is an equivalence class in~$\mathcal E(\,\overline{\mathcal
A}\,)$, the two sets must actually be equal. The proof is
complete.\medskip

{\bf Corollary 7.} {\it If $\mathcal A=\mathcal K$ is compact and
${\rm cap}\,(\mathcal K,a,g)<\infty$, then any solution to the
$\mathcal E(\mathcal K,a\,{\rm cap}\,\mathcal K,g)$-problem gives,
as well, a solution to the $\Gamma(\mathcal
K,a,g)$-problem.}\medskip

{\bf Proof.} This follows from~(\ref{desc2}), when combined
with~(\ref{MMprime}) and~(\ref{S}) for~$a\,{\rm cap}\,\mathcal K$ in
place of~$a$.\medskip

{\bf Remark 12.} Assume ${\rm cap}\,\mathcal A<\infty$, and fix
$\omega\in\mathcal G(\mathcal A,a,g)$ and
$\hat{\omega}\in\hat{\mathcal G}(\mathcal A,a,g)$. Since,
by~(\ref{desc}) and~(\ref{desc2}),
$\kappa(x,\omega)=\kappa(x,\hat\omega)$ nearly everywhere
in~$\mathbf X$, the numbers~$c_i(\omega)$, $i\in I$,
satisfying~(\ref{adm1}) and~(\ref{adm2}) for~$\nu=\omega$, are
actually equal to~$c_i(\hat\omega)$. This implies that
relations~\mbox{(\ref{snt1})\,--\,(\ref{1''})} do hold, as well,
for~$\omega$ in place of~$\hat\omega$.\medskip

{\bf Remark 13.} Observe that, in all the preceding assertions, we
have not imposed any restrictions on the topology of~$A_i$, $i\in
I$. So, all the $\hat{\Gamma}(\mathcal A,a,g)$-,
$\hat{\Gamma}_*(\mathcal A,a,g)$-, and $\Gamma(\mathcal
A,a,g)$-problems are solvable even for a nonclosed
condenser~$\mathcal A$.\medskip

{\bf Remark 14.} If $I=\{1\}$ and $g=1$, Theorems~2\,--\,4 and
Corollary~6 can be derived from~\cite{F1}. Moreover, then one can
choose $\gamma\in\mathcal G(\mathcal A,a,g)$ so that
$$\gamma(\mathbf X)=a_1C(A_1),$$ and exactly this kind of measures
was called by B.~Fuglede {\it interior capacitary distributions
associated with the set\/}~$A_1$. However, this fact in general can
not be extended to the case $I\ne\{1\}$; that is, in general,
$$
\mathcal G(\mathcal A,a,g)\cap\mathcal E(\,\overline{\mathcal
A},a\,{\rm cap}\,\mathcal A,g)=\varnothing,
$$
which can be seen from the unsolvability of the $\mathcal E(\mathcal
A,a\,{\rm cap}\,\mathcal A,g)$-problem.\bigskip

{\bf\Large 7. Interior capacitary constants associated with a
condenser}\nopagebreak\medskip

{\bf 7.1.} Throughout Sec.~7, it is always required that ${\rm
cap}\,(\mathcal A,a,g)<\infty$. Due to the uniqueness statement in
Theorem~3, the following notion naturally arises.\medskip

{\bf Definition 5.} The numbers
$$C_i:=C_i(\mathcal A,a,g):=c_i(\hat\omega),\quad i\in I,$$
satisfying both the relations~(\ref{adm1}) and~(\ref{adm2}) for
$\hat{\omega}\in\hat{\mathcal G}(\mathcal A,a,g)$, are said to be
the ({\it interior}) {\it capacitary constants\/} associated with an
$(m,p)$-condenser~$\mathcal A$.\medskip

{\bf Corollary 8.} {\it The interior capacity ${\rm cap}\,(\mathcal
A,a,g)$ equals the infimum of~$\kappa(\nu,\nu)$, where $\nu$ ranges
over the class of all $\nu\in\mathcal E$} ({\it similarly,}
$\nu\in\mathcal E(\,\overline{\mathcal A}\,)$) {\it such that
$$
\alpha_ia_i\kappa(x,\nu)\geqslant C_i(\mathcal
A,a,g)\,g(x)\quad\mbox{n.\,e.~in \ } A_i,\quad i\in I.
$$
The infimum is attained at any} $\hat\omega\in\hat{\mathcal
G}(\mathcal A,a,g)$ ({\it respectively,} $\omega\in\mathcal
G(\mathcal A,a,g)$), {\it and hence it is an actual
minimum.}\medskip

{\bf Proof.} This follows immediately from Theorems~2\,--\,4 and
Remark~12.\medskip

{\bf 7.2.} Some properties of the interior capacitary constants
$C_i(\mathcal A,a,g)$, $i\in I$, have already been provided by
Theorem~3 and Corollaries~4,~5. Also observe that, if $I$ is a
singleton, then certainly $C_1(\mathcal A,a,g)=1$
(cf.~\cite[Th.~4.1]{F1}).\medskip

{\bf Corollary 9.} {\it $C_i(\,\cdot\,,a,g)$, $i\in I$, are
continuous under exhaustion of~$\mathcal A$ by the increasing family
of all compact condensers $\mathcal K\prec\mathcal A$. Namely,}
$$
C_i(\mathcal A,a,g)=\lim_{\mathcal K\uparrow\mathcal
A}\,C_i(\mathcal K,a,g).$$

{\bf Proof.} Under our assumptions, $0<{\rm cap}\,\mathcal K<\infty$
for every $\mathcal K\in\{\mathcal K\}_\mathcal A$, and hence there
exists $\lambda_\mathcal K\in\mathcal S(\mathcal K,a\,{\rm
cap}\,\mathcal K,g)$. Substituting it into~(\ref{snt1}) yields
\begin{equation}
C_i(\mathcal K,a,g)=\alpha_i\,{\rm cap}\,\mathcal
K^{-1}\,\kappa(\lambda^i_\mathcal K,\lambda_\mathcal K),\quad i\in
I.\label{cntt}
\end{equation}
On the other hand, by Lemma~2 the net
$$
{\rm cap}\,\mathcal A\,{\rm cap}\,\mathcal K^{-1}\,\lambda_\mathcal
K,\quad\mbox{where \ } \mathcal K\in\{\mathcal K\}_\mathcal A,
$$
belongs to the class $\mathbb M(\mathcal A,a\,{\rm cap}\,\mathcal
A,g)$. Substituting it into~(\ref{cnt}) and then combining the
relation obtained with~(\ref{cntt}), we get the corollary.\medskip

In the following assertion we suppose $g_{\min}>0$. According to our
agreement (see~Sec.~4.6), this does hold automatically whenever
$I^-$ is nonempty.\medskip

{\bf Corollary 10.} {\it Assume $C(A_j)=\infty$ for some $j\in I$.
Then
\begin{equation}
C_j(\mathcal A,a,g)\leqslant0.\label{cor4}
\end{equation}
Hence, $C_j(\mathcal A,a,g)=0$ if moreover
$I^-=\varnothing$.}\medskip

{\bf Proof.} Assume, on the contrary, that $C_j>0$. Given
$\hat{\omega}\in\hat{\mathcal G}(\mathcal A,a,g)$, then
$$
\alpha_j a_j\kappa(x,\hat{\omega})\geqslant C_j\,g_{\rm
min}>0\quad\mbox{n.\,e. in \ } A_j,
$$
and therefore, by \cite[Lemma 3.2.2]{F1},
$$
C(A_j)\leqslant a_j^2\,\|\hat{\omega}\|^2\,C_j^{-2}\,g_{\rm
min}^{-2}<\infty,
$$
which is a contradiction. What is left is to show that
$C_j\geqslant0$ provided $I^-=\varnothing$. But this is obvious
because of~(\ref{snt1}).\medskip

{\bf Remark 15.} Observe that the necessity part of Lemma~4, which
has been proved above with elementary arguments, can also be
obtained as a consequence of Corollary~10. Indeed, if (\ref{lemma6})
were true, then by~(\ref{cor4}) the sum of $C_i$, where $i$ ranges
over~$I$, would be not greater than~$0$, which is
impossible.\bigskip

{\bf\Large 8. Interior capacitary distributions associated with a
condenser}\nopagebreak\medskip

As always, we are keeping all our standing assumptions, stated
in~Sec.~4.6. Throughout Sec.~8, it is also required that ${\rm
cap}\,\mathcal A<\infty$.

Our next purpose is to introduce a notion of interior capacitary
distributions~$\gamma_\mathcal A$ associated with a
condenser~$\mathcal A$ such that the distributions obtained possess
properties similar to those of interior capacitary distributions
associated with a set. Fuglede's theory of interior capacities of
sets~\cite{F1} serves here as a model case.

{\bf 8.1.} If $\mathcal A=\mathcal K$ is compact, then, as follows
from Theorem~4, Corollary~7 and Remark~12, any
minimizer~$\lambda_\mathcal K$ in the $\mathcal E(\mathcal K,a\,{\rm
cap}\,\mathcal K,g)$-problem has the desired properties, and so
$\gamma_\mathcal K$ might be defined as $$\gamma_\mathcal
K:=\lambda_\mathcal K,\quad\mbox{where}\quad\lambda_\mathcal
K\in\mathcal S(\mathcal K,a\,{\rm cap}\,\mathcal K,g).$$ However, as
is seen from Remark~14, in the noncompact case the desired notion
can not be obtained as just a direct generalization of the
corresponding one from the theory of interior capacities of sets.
Having in mind that, similar to our model case, the required
distributions should give a solution to the $\Gamma(\mathcal
A,a,g)$-problem and be strongly and $\mathcal A$-vaguely continuous
under exhaustion of~$\mathcal A$ by compact condensers, we arrive at
the following definition.\medskip

{\bf Definition 6.} We shall call $\gamma_\mathcal A\in\mathcal
E(\,\overline{\mathcal A}\,)$ an ({\it interior}) {\it capacitary
distribution\/} associated with~$\mathcal A$ if there exists a
subnet $(\mathcal K_s)_{s\in S}$ of $(\mathcal K)_{\mathcal
K\in\{\mathcal K\}_\mathcal A}$ and
$$\lambda_{\mathcal K_s}\in\mathcal S(\mathcal K_s,a\,{\rm
cap}\,\mathcal K_s,g),\quad s\in S,$$ such that $(\lambda_{\mathcal
K_s})_{s\in S}$ converges to~$\gamma_\mathcal A$ in both the
$\mathcal A$-vague and the strong topologies. Let $\mathcal
D(\mathcal A,a,g)$ denote the collection of all
those~$\gamma_\mathcal A$.\medskip

Application of Lemmas 2 and 5 enables us to rewrite the above
definition in the following, apparently weaker, form:
\begin{equation}
\mathcal D(\mathcal A,a,g)=\mathcal M_0(\mathcal A,a\,{\rm
cap}\,\mathcal A,g).\label{D}
\end{equation}

{\bf Theorem 5.} {\it $\mathcal D(\mathcal A,a,g)$ is nonempty,
contained in an equivalence class in~$\mathcal
E(\,\overline{\mathcal A}\,)$, and compact in the induced $\mathcal
A$-vague topology. Furthermore,
\begin{equation}
\mathcal D(\mathcal A,a,g)\subset\mathcal G(\mathcal
A,a,g)\cap\mathcal E(\,\overline{\mathcal A},\leqslant\!a\,{\rm
cap}\,\mathcal A,g).\label{gammain}
\end{equation}
Given its element $\gamma:=\gamma_\mathcal A$, then
\begin{equation}
\|\gamma\|^2={\rm cap}\,\mathcal A, \label{5.1}
\end{equation}
\begin{equation}
\alpha_i a_i\kappa(x,\gamma)\geqslant C_i\,g(x)\quad\mbox{n.\,e. in
\ } A_i,\quad i\in I,\label{5.3}
\end{equation}
where $C_i=C_i(\mathcal A,a,g)$, $i\in I$, are the interior
capacitary constants. Actually,
\begin{equation}
C_i=\frac{\alpha_i\kappa(\gamma^i,\gamma)}{{\rm cap}\,\mathcal
A}="\!\inf_{x\in A_i}\!"\,\,\frac{\alpha_i
a_i\kappa(x,\gamma)}{g(x)}\,,\quad i\in I.\label{5.2}
\end{equation}
If $I^-\ne\varnothing$, assume moreover that the kernel
$\kappa(x,y)$ is continuous for $x\ne y$, while\/}
$\kappa(\cdot,y)\to0$ ({\it as\/}~$y\to\infty$) {\it uniformly on
compact sets. Then, for every $i\in I$,}
\begin{align}
\alpha_i a_i\kappa(x,\gamma)&\leqslant
C_i\,g(x)&&\!\!\!\!\!\!\!\!\!\!\!\!\!\!\!\!\!\!\!\!\!\!\!\!\!\!\!\!\!\!\text{\it
for all \ } x\in
S(\gamma^i),\label{5.4}\\
\intertext{\it and hence} \alpha_i
a_i\kappa(x,\gamma)&=C_i\,g(x)&&\!\!\!\!\!\!\!\!\!\!\!\!\!\!\!\!\!\!\!\!\!\!\!\!\!\!\!\!\!\!\text{\it
n.\,e. in \ } A_i\cap S(\gamma^i).\notag
\end{align}

Thus, an interior capacitary distribution $\gamma_\mathcal A$ is
{\it unique\/} if the kernel is additionally assumed to be strictly
positive definite and all $\overline{A_i}$, $i\in I$, are mutually
disjoint.\medskip

{\bf Remark 16.} As is seen from the preceding theorem, the
properties of interior capacitary distributions associated with a
condenser are quite similar to those of interior capacitary
distributions associated with a set (cf.~\cite[Th.~4.1]{F1}). The
only important difference is that the sign~$\leqslant$ in the
inclusion~(\ref{gammain}) in general can not be omitted~--- even for
a closed, noncompact condenser.\medskip

{\bf Remark 17.} Like as in the theory of interior capacities of
sets, in general none of the $i$-coordinates of $\gamma_\mathcal A$
is concentrated on~$A_i$ (unless $A_i$ is closed). Indeed, let
$\mathbf X=\mathbb R^n$, $n\geqslant3$, $\kappa(x,y)=|x-y|^{2-n}$,
$g=1$, $I^+=\{1\}$, $I^-=\{2\}$, $a_1=a_2=1$, and let $A_1=\{x:
|x|<r\}$ and $A_2=\{x: |x|>R\}$, where $0<r<R<\infty$. Then it can
be shown that
$$\gamma_\mathcal
A=\gamma_{\overline{\mathcal A}}=\bigl[\theta^+-\theta^-\bigr]\,{\rm
cap}\,\mathcal A,$$ where $\theta^+$ and~$\theta^-$ are obtained by
the uniform distribution of unit mass over the spheres~$S(0,r)$
and~$S(0,R)$, respectively. Hence, $|\gamma_\mathcal
A|(A)=0$.\medskip

{\bf 8.2.} The purpose of this section is to point out
characteristic properties of the interior capacitary distributions
and the interior capacitary constants.\medskip

{\bf Proposition 1.} {\it Assume $\mu\in\mathcal
E(\,\overline{\mathcal A}\,)$ has the properties
$$
\|\mu\|^2={\rm cap}\,(\mathcal A,a,g),
$$
$$
\alpha_i
a_i\kappa(x,\mu)\geqslant\frac{\alpha_i\kappa(\mu^i,\mu)}{{\rm
cap}\,\mathcal A}\,g(x)\quad\mbox{\it n.\,e. in \ } A_i,\quad i\in
I.
$$
Then $\mu$ is equivalent in $\mathcal E(\,\overline{\mathcal A}\,)$
to every $\gamma_\mathcal A\in\mathcal D(\mathcal A,a,g)$, and for
all $i\in I$,}
$$
C_i(\mathcal A,a,g)=\frac{\alpha_i\kappa(\mu^i,\mu)}{{\rm
cap}\,\mathcal A}="\!\inf_{x\in A_i}\!"\,\,\frac{\alpha_i
a_i\kappa(x,\mu)}{g(x)}\,.
$$

Actually, there holds the following stronger result, to be proved in
Sec.~16 below.\medskip

{\bf Proposition 2.} {\it Let $\nu\in\mathcal E(\,\overline{\mathcal
A}\,)$ and $\tau_i\in\mathbb R$, $i\in I$, satisfy the relations
\begin{equation}
\alpha_i a_i\kappa(x,\nu)\geqslant\tau_i\,g(x)\quad\mbox{\it n.\,e.
in \ } A_i,\quad i\in I,\label{desc3'}
\end{equation}
\begin{equation}
\sum_{i\in I}\,\tau_i=\frac{{\rm cap}\,\mathcal A+\|\nu\|^2}{2\,{\rm
cap}\,\mathcal A}\,.\label{adm2''}
\end{equation}
Then $\nu$ is equivalent in $\mathcal E(\,\overline{\mathcal A}\,)$
to every $\gamma_\mathcal A\in\mathcal D(\mathcal A,a,g)$, and for
all $i\in I$,}
\begin{equation}
\tau_i=C_i(\mathcal A,a,g)="\!\inf_{x\in A_i}\!"\,\,\frac{\alpha_i
a_i\kappa(x,\nu)}{g(x)}\,.\label{un}
\end{equation}

Thus, under the conditions of Proposition~1 or~2, if moreover
$\kappa$ is strictly positive definite and all $\overline{A_i}$,
$i\in I$, are mutually disjoint, then the measure under
consideration is actually the (unique) interior capacitary
distribution~$\gamma_\mathcal A$.\bigskip

{\bf\Large 9. On continuity of the capacities, capacitary
distributions, and capacitary constants}\nopagebreak\medskip

{\bf 9.1.} Given $\mathcal A_n=(A_i^n)_{i\in I}$, $n\in\mathbb N$,
and $\mathcal A$ in $\mathfrak C_{m,p}$, write $\mathcal
A_n\uparrow\mathcal A$ if $\mathcal A_n\prec\mathcal A_{n+1}$ for
all~$n$ and
$$A_i=\bigcup_{n\in\mathbb N}\,A_i^n,\quad i\in I.$$

Following \cite[Chap.~1, \S\,9]{B1}, we call a locally compact space
{\it countable at infinity\/} if it can be written as a countable
union of compact sets.\medskip

{\bf Theorem 6.} {\it Suppose that either $g_{\rm min}>0$ or the
space $\mathbf X$ is countable at infinity. If $\mathcal A_n$,
$n\in\mathbb N$, are universally measurable and $\mathcal
A_n\uparrow\mathcal A$, then
\begin{equation}
{\rm cap}\,(\mathcal A,a,g)=\lim_{n\in\mathbb N}\,{\rm
cap}\,(\mathcal A_n,a,g). \label{6.1}
\end{equation}
Assume moreover ${\rm cap}\,(\mathcal A,a,g)$ to be finite, and let
$\gamma_n:=\gamma_{\mathcal A_n}$, $n\in\mathbb N$, denote an
interior capacitary distribution associated with~$\mathcal A_n$. If
$\gamma$ is an $\mathcal A$-vague limit point
of\/}~$(\gamma_n)_{n\in\mathbb N}$ ({\it such a $\gamma$ exists\/}),
{\it then $\gamma$ is actually an interior capacitary distribution
associated with the condenser~$\mathcal A$, and
$$ \lim_{n\in\mathbb
N}\,\|\gamma_n-\gamma\|^2=0.
$$
Furthermore,}
\begin{equation}
C_i(\mathcal A,a,g)=\lim_{n\in\mathbb N}\,C_i(\mathcal
A_n,a,g),\quad i\in I\smallskip.\label{6.3}
\end{equation}

Thus, if $\kappa$ is additionally assumed to be strictly positive
definite (hence, perfect) and all $\overline{A_i}$, $i\in I$, are
mutually disjoint, then the (unique) interior capacitary
distribution associated with~$\mathcal A_n$ converges both $\mathcal
A$-vaguely and strongly to the (unique) interior capacitary
distribution associated with~$\mathcal A$.\medskip

{\bf Remark 18.} Theorem 6 remains true if $(\mathcal
A_n)_{n\in\mathbb N}$ is replaced by the increasing ordered family
of all compact condensers~$\mathcal K$ such that $\mathcal
K\prec\mathcal A$. Moreover, then the assumption that either $g_{\rm
min}>0$ or $\mathbf X$ is countable at infinity can be omitted. Cf.,
e.\,g., Lemma~2 and Corollary~9.\medskip

{\bf Remark 19.} If $I=\{1\}$ and $g=1$, Theorem~6 has been proved
in~\cite[Th.~4.2]{F1}.\medskip

{\bf 9.2.} The remainder of the article is devoted to proving the
results formulated in Sec.~\mbox{5\,--\,9} and is organized as
follows. Theorems~2,~3,~5, and~6 are proved in~Sec.~14,~15, and~17.
Their proofs utilize the description of the potentials of measures
from the classes $\mathcal M'(\mathcal A,a,g)$ and $\mathcal
M_0(\mathcal A,a,g)$, to be given in Sec.~12 and~13 by Lemmas~9
and~10. In turn, Lemmas~9 and~10 use the theorem on the strong
completeness of proper subspaces of~$\mathcal E$, which is a subject
of Sec.~10.\bigskip

{\bf\Large 10. On the strong completeness}\nopagebreak\medskip

{\bf 10.1.} Keeping all our standing assumptions on~$\kappa$, $g$,
and~$\mathcal A$, stated in~Sec.~4.6, we consider $\mathcal
E(\,\overline{\mathcal A},\leqslant\!a,g)$ to be a topological
subspace of the semimetric spa\-ce~$\mathcal E(\,\overline{\mathcal
A}\,)$; the induced topology is likewise called the {\it strong\/}
topology.\medskip

{\bf Theorem 7.} {\it Suppose $\mathcal A$ is closed. Then the
semimetric space $\mathcal E(\mathcal A,\leqslant\!a,g)$ is
complete. In more detail, if $(\mu_s)_{s\in S}\subset\mathcal
E(\mathcal A,\leqslant\!a,g)$ is a strong Cauchy net and $\mu$~is
its $\mathcal A$-vague cluster point\/} ({\it such a $\mu$
exists\/}), {\it then $\mu\in\mathcal E(\mathcal A,\leqslant\!a,g)$
and
\begin{equation}
\lim_{s\in S}\,\|\mu_s-\mu\|^2=0.\label{str}
\end{equation}
Assume, in addition, that the kernel is strictly positive definite
and all $A_i$, $i\in I$, are mutually disjoint. If moreover
$(\mu_s)_{s\in S}\subset\mathcal E(\mathcal A,\leqslant\!a,g)$
converges strongly to $\mu_0\in\mathcal E(\mathcal A)$, then
actually $\mu_0\in\mathcal E(\mathcal A,\leqslant\!a,g)$ and
$\mu_s\to\mu_0$ $\mathcal A$-vaguely}.\medskip

{\bf Remark 20.} This theorem is certainly of independent interest
since, according to the well-known counterexample by
H.~Cartan~\cite{Car}, the pre-Hilbert space~$\mathcal E$ is strongly
incomplete even for the Newton kernel $|x-y|^{2-n}$ in~$\mathbb
R^n$, $n\geqslant3$.\medskip

{\bf Remark 21.} Assume the kernel is strictly positive definite
(hence, perfect). If moreover $I^-=\varnothing$, then Theorem~7
remains valid for $\mathcal E(\mathcal A)$ in place of $\mathcal
E(\mathcal A,\leqslant\!a,g)$ (cf.~Theorem~1).  A question still
unanswered is {\bf whether this is the case if $I^+$ and $I^-$ are
both nonempty\/}. We can however show that this is really so for the
Riesz kernels $|x-y|^{\alpha-n}$, $0<\alpha<n$, in~$\mathbb R^n$,
$n\geqslant2$ (cf.~\cite[Th.~1]{Z2}). The proof utilizes Deny's
theorem~\cite{D1} stating that, for the Riesz kernels, $\mathcal E$
can be completed with making use of distributions of finite
energy.\medskip

{\bf 10.2.} We start by auxiliary assertions to be used in the proof
of Theorem~7.\medskip

{\bf Lemma 7.} \ {\it $\mathcal E(\mathcal A,\leqslant\!a,g)$ is
$\mathcal A$-vaguely bounded.}\medskip

{\bf Proof.} Let $K\subset A_i$, $i\in I$, be compact. Since $g$ is
positive and continuous, the inequalities
$$
a_i\geqslant\int g\,d\mu^i\geqslant\mu^i(K)\,\min_{x\in K}\
g(x),\quad\mbox{where \ }\mu\in\mathcal E(\mathcal
A,\leqslant\!a,g),
$$
yield
$$
\sup_{\mu\in\mathcal E(\mathcal A,\leqslant a,g)}\,\mu^i(K)<\infty,
$$
and the lemma follows.\medskip

{\bf Lemma 8.} {\it Suppose $\mathcal A$ is closed. If a net
$(\mu_s)_{s\in S}\subset\mathcal E(\mathcal A,\leqslant\!a,g)$ is
strongly bounded, then its $\mathcal A$-vague cluster set is
nonempty and contained in $\mathcal E(\mathcal
A,\leqslant\!a,g)$.}\medskip

{\bf Proof.} We begin by showing that the nets $(\mu_s^i)_{s\in S}$,
$i\in I$, are strongly bounded as well, i.\,e.,
\begin{equation}
\sup_{s\in S}\ \|\mu_s^i\|<\infty,\quad
i\in I.  \label{7.1}
\end{equation}
This is obvious when $I^-=\varnothing$ and $\kappa\geqslant0$;
hence, one can assume that either $I^-\ne\varnothing$, or
$I^-=\varnothing$ while $\mathbf X$ is compact. In any case,
both~(\ref{g}) and~(\ref{bou}) hold. Since
\begin{equation}
\int g\,d\mu_s^i\leqslant a_i,\quad i\in I,
\label{oh}
\end{equation}
(\ref{g}) implies
\begin{equation} \sup_{s\in S}\
\mu_s^i(\mathbf X)\leqslant a_i\,g_{\min}^{-1}<\infty,\quad i\in I.
\label{7.3}
\end{equation}
When combined with~(\ref{bou}), this shows that
$\kappa(\mu^+_s,\mu^-_s)$ remains bounded from above on~$S$, and
hence so do $\|\mu^+_s\|^2$ and $\|\mu^-_s\|^2$. Since $\kappa$ is
bounded from below on $\mathbf X\times\mathbf X$, repeated
application of~(\ref{7.3}) gives~(\ref{7.1}) as desired.

Moreover, for every $i\in I$, $(\mu_s^i)_{s\in S}$ is vaguely
bounded according to the preceding lemma, while $\mathfrak M^+(A_i)$
is vaguely closed. Since any vaguely bounded part of~$\mathfrak M$
is vaguely relatively compact (see, e.\,g., \cite[Chap.~III, \S~2,
Prop.~9]{B2}), there exists a vague cluster point of
$(\mu^i_s)_{s\in S}$, say~$\mu^i$, and $\mu^i\in\mathfrak M^+(A_i)$.

It remains to show that $\mu^i$ is of finite energy and
satisfies~(\ref{oh}) with~$\mu^i$ in place of~$\mu_s^i$. To this
end, recall that, if $\mathbf Y$~is a locally compact Hausdorff
space and $\psi$~is a lower semicontinuous function on~$\mathbf Y$
such that $\psi\geqslant0$ (unless its support is compact), then the
map
$$
\nu\mapsto\int\psi\,d\nu,\quad\nu\in\mathfrak M^+(\mathbf Y),
$$
is lower semicontinuous in the induced vague topology (see,
e.\,g.,~\cite{F1}). Applying this to $\mathbf Y=A_i\times A_i$,
$\psi=\kappa|_{A_i\times A_i}$ and, subsequently, $\mathbf Y=A_i$,
$\psi=g|_{A_i}$, we derive the required properties of $\mu^i$
from~(\ref{7.1}) and~(\ref{oh}).\medskip

{\bf 10.3. Proof of Theorem 7.} Suppose $\mathcal A$ is closed, and
let $(\mu_s)_{s\in S}$ be a strong Cauchy net in $\mathcal
E(\mathcal A,\leqslant\!a,g)$. Since such a net converges strongly
to every its strong cluster point, $(\mu_s)_{s\in S}$ can certainly
be assumed to be strongly bounded. Then, by Lemma~8, there exists an
$\mathcal A$-vague cluster point~$\mu$ of~$(\mu_s)_{s\in S}$, and
\begin{equation}
\mu\in\mathcal E(\mathcal A,\leqslant\!a,g).
\label{leqslant}
\end{equation}
We next proceed
to verify (\ref{str}).

Without loss of generality we can also assume that, for every $i\in
I$,
$$
\mu^i_s\to\mu^i\quad\mbox{vaguely}.
$$ Since, by~(\ref{7.1}),
$(\mu^i_s)_{s\in S}$ is strongly bounded, the property~$(CW)$
(see~Sec.~2) shows that $\mu^i_s$ approaches~$\mu^i$ in the weak
topology as well, and so
$$
R\mu_s\to R\mu\quad\mbox{weakly}.
$$
This gives
$$
\|\mu_s-\mu\|^2=\|R\mu_s-R\mu\|^2=\lim_{l\in
S}\,\kappa(R\mu_s-R\mu,R\mu_s-R\mu_l),
$$
and hence, by the Cauchy-Schwarz inequality,
$$\|\mu_s-\mu\|^2\leqslant
\|\mu_s-\mu\|\,\liminf_{l\in S}\,\|\mu_s-\mu_l\|,
$$
which proves (\ref{str}) as required, because $\|\mu_s-\mu_l\|$
becomes arbitrarily small when $s,\,l\in S$ are both large enough.

Suppose now that $\kappa$ is strictly positive definite, while all
$A_i$, $i\in I$, are mutually disjoint, and let the net
$(\mu_s)_{s\in S}$ converge strongly to some $\mu_0\in\mathcal
E(\mathcal A)$. Given an vague limit point~$\mu$ of~$(\mu_s)_{s\in
S}$, then we conclude from~(\ref{str}) that $\|\mu_0-\mu\|=0$, hence
$\mu_0\cong\mu$ since $\kappa$ is strictly positive definite, and
finally $\mu_0\equiv\mu$ because $A_i$, $i\in I$, are mutually
disjoint. In view of~(\ref{leqslant}), this means that
$\mu_0\in\mathcal E(\mathcal A,\leqslant\!a,g)$, which is a part of
the desired conclusion.

Moreover, $\mu_0$ has thus been shown to be identical to any
$\mathcal A$-vague cluster point of~$(\mu_s)_{s\in S}$. Since the
vague topology is separated, this implies that $\mu_0$ is actually
its $\mathcal A$-vague limit (cf.~\cite[Chap.~I, \S~9, n$^\circ$\,1,
cor.]{B1}), which completes the proof.\bigskip

{\bf\Large 11. Proof of Lemma 5}\nopagebreak\medskip

Fix any $(\mu_s)_{s\in S}$ and $(\nu_t)_{t\in T}$ in $\mathbb
M(\mathcal A,a,g)$. It follows by standard arguments that
\begin{equation}
\lim_{(s,t)\in S\times T}\,\|\mu_s-\nu_t\|^2=0, \label{fund}
\end{equation}
where $S\times T$ denotes the directed product of the directed
sets~$S$ and~$T$ (see, e.\,g.,~\cite[Chap.~2,~\S~3]{K}). Indeed, by
the convexity of the class $\mathcal E(\mathcal A,a,g)$,
$$
2\,\|\mathcal E(\mathcal A,a,g)\|
\leqslant{\|\mu_s+\nu_t\|}\leqslant\|\mu_s\|+\|\nu_t\|,
$$
and hence, by (\ref{min}),
$$
\lim_{(s,t)\in S\times T}\,\|\mu_s+\nu_t\|^2=4\,\|\mathcal
E(\mathcal A,a,g)\|^2.
$$
Then the parallelogram identity gives~(\ref{fund}) as claimed.

Relation~(\ref{fund}) implies that $(\mu_s)_{s\in S}$ is strongly
fundamental. Therefore Theorem~7 shows that there exists an
$\mathcal A$-vague cluster point~$\mu_0$ of~$(\mu_s)_{s\in S}$, and
moreover $\mu_0\in\mathcal E(\,\overline{\mathcal
A},\leqslant\!a,g)$ and $\mu_s\to\mu_0$ strongly. This means that
$\mathcal M(\mathcal A,a,g)$ and $\mathcal M'(\mathcal A,a,g)$ are
both nonempty and satisfy the inclusion~(\ref{MMprime}).

What is left is to prove that $\mu_s\to\chi$ strongly, where
$\chi\in\mathcal M'(\mathcal A,a,g)$ is arbitrarily given. But then
one can choose a net in~$\mathbb M(\mathcal A,a,g)$ converging
to~$\chi$ strongly, and repeated application of~(\ref{fund}) leads
immediately to the desired conclusion.\bigskip

{\bf\Large 12. Potentials of strong cluster points of minimizing
nets}\nopagebreak\medskip

{\bf 12.1.} The aim of this section is to provide a description of
the potentials of measures from the class $\mathcal M'(\mathcal
A,a,g)$. As usual, we are keeping all our standing assumptions,
stated in Sec.~4.6.\medskip

{\bf Lemma 9.} {\it There exist $\eta_i\in\mathbb R$, $i\in I$, such
that, for every $\chi\in\mathcal M'(\mathcal A,a,g)$,
\begin{equation}
\alpha_ia_i\kappa(x,\chi) \geqslant\alpha_i\eta_i
g(x)\quad\mbox{n.\,e. in \ } A_i,\quad i\in I, \label{1.1}
\end{equation}
\begin{equation}
\sum_{i\in I}\alpha_i\eta_i=\|\mathcal E(\mathcal A,a,g)\|^2.
\label{1.2}
\end{equation}
These $\eta_i$, $i\in I$, are determined uniquely and given by
either of the formulas
\begin{align}
\eta_i&=\kappa(\zeta^i,\zeta),\label{1.3'}\\[2pt]
\eta_i&=\lim_{s\in S}\,\kappa(\mu_s^i,\mu_s),\label{1.3''}
\end{align}
where $\zeta\in\mathcal M(\mathcal A,a,g)$ and $(\mu_s)_{s\in
S}\in\mathbb M(\mathcal A,a,g)$ are arbitrarily chosen}.\medskip

{\bf Proof.} Throughout the proof, we shall assume every net
$(\mu_s)_{s\in S}\in\mathbb M(\mathcal A,a,g)$ to be strongly
bounded, which certainly involves no loss of generality. Then all
the nets $(\mu^i_s)_{s\in S}$, $i\in I$, are strongly bounded as
well (see~the proof of Lemma~8).

Choose $(\mu_t)_{t\in T}\in\mathbb M(\mathcal A,a,g)$ with the
property that, for every $i\in I$, there exists the limit (finite or
infinite)
\begin{equation}
\eta_i:=\lim_{t\in T}\,\kappa(\mu_t^i,\mu_t). \label{2.1}
\end{equation}
We proceed to show that $\eta_i$, $i\in I$, so defined, satisfy
both~(\ref{1.1}) and~(\ref{1.2}).

Given $\chi\in\mathcal M'(\mathcal A,a,g)$, suppose, contrary to our
claim, that for some $j\in I$ there exists a set $E_j\subset A_j$ of
interior capacity nonzero such that
\begin{equation}
\alpha_ja_j\kappa(x,\chi)<\alpha_j\eta_jg(x)\quad\mbox{for all \ }
x\in E_j. \label{3.2} \end{equation} Then one can choose
$\nu\in\mathcal E^+$ with compact support so that $S(\nu)\subset
E_j$ and
$$\int g\,d\nu=a_j.$$
Integrating the inequality in (\ref{3.2}) with respect to~$\nu$
gives
\begin{equation}
\alpha_j\,\bigl[\kappa(\chi,\nu)-\eta_j\bigr]<0.
\label{3.3} \end{equation}

To get a contradiction, for every $\tau\in(0,1]$ write
$$
\tilde{\mu}^i_t:=\left\{
\begin{array}{cl} \mu^j_t-\tau\bigl(\mu_t^j-\nu\bigr) & \mbox{if \ }
i=j,\\[4pt]
\mu^i_t & \mbox{otherwise}.\\ \end{array} \right.
$$
Clearly,
$$
\tilde{\mu}_t:=\sum_{i\in I}\,\alpha_i\tilde{\mu}^i_t\in\mathcal
E^0(\mathcal A,a,g),\quad t\in T,
$$
and consequently
\begin{equation} \|\mathcal
E(\mathcal A,a,g)\|^2\leqslant \|\tilde{\mu}_t\|^2=\|\mu_t\|^2
-2\alpha_j\tau\,\kappa(\mu_t,\mu_t^j-\nu)+\tau^2\|\mu_t^j-\nu\|^2.
\label{3.4}
\end{equation}
The coefficient of~$\tau^2$ is bounded from above on~$T$ (say
by~$M_0$), while by Lemma~5
$$
\lim_{t\in T}\,\|\mu_t-\chi\|^2=0.
$$
From (\ref{2.1}) and~(\ref{3.4}) we therefore obtain
$$
0\leqslant
M_0\tau^2+2\alpha_j\tau\,\bigl[\kappa(\chi,\nu)-\eta_j\bigr].
$$
By letting here $\tau$ tend to~$0$, we arrive at a contradiction to
(\ref{3.3}).

It has thus been proved that $\eta_i$, $i\in I$, defined by means
of~(\ref{2.1}), satisfy~(\ref{1.1}). Note that
$\kappa(\,\cdot\,,R\chi)$, being the potential of a measure of
finite energy, is finite nearly everywhere in~$\mathbf X$
(see~\cite{F1}), and hence so is $\kappa(\,\cdot\,,\chi)$. Since, by
Lemma~3, $C(A_i)>0$ for all $i\in I$, it follows from~(\ref{1.1})
that
$$\alpha_i\eta_i<\infty,\quad i\in I.$$ Hence, $\sum_{i\in
I}\alpha_i\eta_i$ is well defined and, by~(\ref{2.1}),
$$
\sum_{i\in I}\alpha_i\eta_i=\lim_{t\in T}\,\|\mu_t\|^2=\|\mathcal
E(\mathcal A,a,g)\|^2.
$$
This means that $\eta_i$, $i\in I$, are finite and satisfy
also~(\ref{1.2}) as required.

To prove the statement on uniqueness, consider some other $\eta'_i$,
$i\in I$, satisfying both~(\ref{1.1}) and~(\ref{1.2}). Then they are
necessarily finite, and for every~$i$,
\begin{equation}
\alpha_ia_i\kappa(x,\chi)\geqslant\max\bigl\{\alpha_i\eta_i,\,
\alpha_i\eta'_i\bigr\}\,g(x)\quad\mbox{n.\,e. in \ } A_i,\label{4.1}
\end{equation}
which follows from the property of subadditivity of $C(\,\cdot\,)$,
mentioned in~Sec.~6.1. Since $\mu_t^i$ is concentrated on~$A_i$ and
has finite energy and compact support, application of~\cite[
Lemma~2.3.1]{F1} shows that the inequality in~(\ref{4.1}) holds
$\mu_t^i$-almost everywhere in~$\mathbf X$. Integrating it with
respect to~$\mu_t^i$ and then summing up over all $i\in I$, in view
of $\int g\,d\mu_t^i=a_i$ we have
$$
\kappa(\mu_t,\chi)\geqslant\sum_{i\in I}\,
\max\bigl\{\alpha_i\eta_i,\,\alpha_i\eta'_i\bigr\},\quad t\in T.
$$
Passing here to the limit as $t$ ranges over $T$, we get
$$
\|\chi\|^2=\lim_{t\in T}\kappa(\mu_t,\chi)\geqslant \sum_{i\in
I}\max\bigl\{\alpha_i\eta_i,\,\alpha_i\eta'_i\bigr\}\geqslant\sum_{i\in
I}\alpha_i\eta_i= \|\mathcal E(\mathcal A,a,g)\|^2,
$$
and hence
$$
\max\bigl\{\alpha_i\eta_i,\,\alpha_i\eta'_i\bigr\}
=\alpha_i\eta_i,\quad i\in I,
$$
for the extreme left and right parts of the above chain of
inequalities are equal. Applying the same arguments again, but with
the roles of $\eta_i$ and $\eta'_i$ reversed, we conclude that
$\eta_i=\eta'_i$ for all $i\in I$, as claimed.

It remains to show that $\eta_i$, $i\in I$, can be written in the
form~(\ref{1.3'}) or~(\ref{1.3''}). To this end, fix $(\mu_s)_{s\in
S}\in\mathbb M(\mathcal A,a,g)$. Then it follows at once from the
above reasoning that, for every $i\in I$, any cluster point of the
net $\kappa(\mu_s^i,\mu_s)$, $s\in S$, coincides with~$\eta_i$.
Hence, there exists $\lim_{s\in S}\,\kappa(\mu_s^i,\mu_s)$ and it
equals~$\eta_i$.

Passing to a subnet if necessary, by Lemma~8 we can also assume
$(\mu_s)_{s\in S}$ to be $\mathcal A$-vaguely convergent, say
to~$\zeta$. The proof will be completed once we prove
\begin{equation}
\kappa(\zeta^i,\zeta)=\lim_{s\in S}\,\kappa(\mu_s^i,\mu_s),\quad
i\in I.\label{proof}
\end{equation}
Since $\|\mu_s^i\|$ is bounded from above on~$S$ (say by~$M_1$),
while $\mu_s^i\to\zeta^i$ vaguely, the property~$(CW)$ yields that
$\mu_s^i$ approaches~$\zeta^i$ also weakly. Hence, for every
$\varepsilon>0$,
$$
|\kappa(\zeta^i-\mu_s^i,\zeta)|<\varepsilon
$$
whenever $s\in S$ is large enough. Furthermore, by the
Cauchy-Schwarz inequality,
$$
|\kappa(\mu_s^i,\zeta-\mu_s)|\leqslant M_1\|\zeta-\mu_s\|,\quad s\in
S.
$$
Since, by Lemma 5, $\mu_s\to\zeta$ strongly, the last two relations
combined give~(\ref{proof}).\medskip

{\bf 12.2.} In what follows, $\eta_i=:\eta_i(\mathcal A,a,g)$, $i\in
I$, will always denote the numbers appeared in~Lemma~9. They are
uniquely determined by relation~(\ref{1.1}), where $\chi\in\mathcal
M'(\mathcal A,a,g)$ is arbitrarily chosen, taken together
with~(\ref{1.2}). This statement on uniqueness can actually be
strengthened as follows.\medskip

{\bf Lemma 9$'$.} {\it Given $\chi\in\mathcal M'(\mathcal A,a,g)$,
choose $\eta'_i$, $i\in I$, so that
$$
\sum_{i\in I}\alpha_i\eta'_i\geqslant\|\mathcal E(\mathcal
A,a,g)\|^2.
$$
If there holds\/}~(\ref{1.1}) {\it for $\eta'_i$ in place
of~$\eta_i$, then $\eta'_i=\eta_i$ for all $i\in I$.}\medskip

{\bf Proof.} This follows in the same manner as the uniqueness
statement in~Lemma~9.\medskip

{\bf 12.3.} The following assertion is specifying Lemma 9 for a
compact condenser~$\mathcal K$.\medskip

{\bf Corollary 11.} {\it Let $\mathcal A=\mathcal K$ be compact.
Given $\lambda_\mathcal K\in\mathcal S(\mathcal K,a,g)$, then for
every~$i$,}
\begin{align}
\alpha_ia_i\kappa(x,\lambda_\mathcal
K)&\geqslant\alpha_i\kappa(\lambda^i_\mathcal K,\lambda_\mathcal
K)\,g(x)&&\!\!\!\!\!\!\!\!\!\!\!\!\!\text{\it n.\,e. in \ } K_i,\label{hh}\\
\intertext{\it and hence} a_i\kappa(x,\lambda_\mathcal K)&
=\kappa(\lambda^i_\mathcal K,\lambda_\mathcal
K)\,g(x)&&\!\!\!\!\!\!\!\!\!\!\!\!\!\text{\it $\lambda^i_\mathcal
K$-almost everywhere.}\label{1.1''}
\end{align}

{\bf Proof.} In view of (\ref{S}) and (\ref{1.3'}), $\eta_i(\mathcal
K,a,g)$, $i\in I$, can be written in the form
$$\eta_i(\mathcal K,a,g)=\kappa(\lambda^i_\mathcal K,\lambda_\mathcal
K),$$ which leads to~(\ref{hh}) when substituted into~(\ref{1.1}).
Since $\lambda^i_\mathcal K$ has finite energy and is supported
by~$K_i$, the inequality in~(\ref{hh}) holds $\lambda^i_\mathcal
K$-almost everywhere in~$\mathbf X$. Hence, (\ref{1.1''}) must be
true, for if not, we would arrive at a contradiction by integrating
the inequality in~(\ref{hh}) with respect to~$\lambda^i_\mathcal
K$.\bigskip

{\bf\Large 13. Potentials of $\mathcal A$-vague cluster points of
minimizing nets}\nopagebreak\medskip

In this section we shall restrict ourselves to measures $\xi$ of the
class $\mathcal M_0(\mathcal A,a,g)$. It is clear from Corollary~3
that their potentials have all the properties described in Lemmas~9
and~9$'$. Our purpose is to show that, under proper additional
restrictions on the kernel, that description can be sharpened as
follows.\medskip

{\bf Lemma 10.} {\it In the case where $I^-\ne\varnothing$, assume
moreover that $\kappa(x,y)$ is continuous for $x\ne y$, while\/}
$\kappa(\cdot,y)\to0$ ({\it as\/}~$y\to\infty$) {\it uniformly on
compact sets. Given $\xi\in\mathcal M_0(\mathcal A,a,g)$, then for
all $i\in I$,}
\begin{align} \alpha_ia_i\kappa(x,\xi)&
\geqslant\alpha_i\kappa(\xi^i,\xi)\,g(x)&&\!\!\!\!\!\!\!\!\!\!\!\!\!\!\!\!\!\!\text{\it
n.\,e. in \ }
A_i,\label{1.11}\\[4pt]
\alpha_i a_i\kappa(x,\xi)&\leqslant
\alpha_i\kappa(\xi^i,\xi)\,g(x)&&\!\!\!\!\!\!\!\!\!\!\!\!\!\!\!\!\!\!\text{\it
for all \ }
x\in S(\xi^i),\label{desc4'}\\
\intertext{\it and hence}
a_i\kappa(x,\xi)&=\kappa(\xi^i,\xi)\,g(x)&&\!\!\!\!\!\!\!\!\!\!\!\!\!\!\!\!\!\!\text{\it
n.\,e. in \ } A_i\cap S(\xi^i).\notag\end{align}

{\bf Proof.} Choose $\lambda_\mathcal K\in\mathcal S(\mathcal
K,a,g)$ such that $\xi$ is an $\mathcal A$-vague cluster point of
the net $(\lambda_\mathcal K)_{\mathcal K\in\{\mathcal K\}_\mathcal
A}$. Since this net belongs to $\mathbb M(\mathcal A,a,g)$, from
(\ref{1.3'}) and~(\ref{1.3''}) we get
$$
\eta_i=\kappa(\xi^i,\xi)=\lim_{\mathcal K\in\{\mathcal K\}_\mathcal
A}\,\kappa(\lambda^i_\mathcal K,\lambda_\mathcal K),\quad i\in I.
$$
Substituting this into~(\ref{1.1}) with~$\xi$ in place of~$\chi$
gives~(\ref{1.11}) as required.

We next proceed to prove~(\ref{desc4'}). To this end, fix $i$ (say
$i\in I^+$) and $x_0\in S(\xi^i)$. Without loss of generality it can
certainly be assumed that
\begin{equation}
\lambda_\mathcal K\to\xi\quad\mbox{$\mathcal A$-vaguely},
\label{l11}\end{equation} since otherwise we shall pass to a subnet
and change the notation. Then, due to~(\ref{1.1''}) and~(\ref{l11}),
there exist $x_\mathcal K\in S(\lambda^i_\mathcal K)$ with the
following properties:
\begin{equation}
x_\mathcal K\to x_0\quad\mbox{as \ } \mathcal K\uparrow\mathcal
A,\label{l112}
\end{equation}
$$
a_i\kappa(x_\mathcal K,\lambda_\mathcal K)=\kappa(\lambda^i_\mathcal
K,\lambda_\mathcal K)\,g(x_\mathcal K).
$$

Taking into account that, by~\cite[Lemma 2.2.1]{F1}, the map
$$(x,\nu)\mapsto\kappa(x,\nu)$$
is lower semicontinuous on $\mathbf X\times\mathfrak M^+$ in the
topology of a Cartesian product (where $\mathfrak M^+$ is equipped
with the vague topology), we conclude from what has already been
shown that the desired relation~(\ref{desc4'}) will follow once we
prove
\begin{equation}
\kappa(x_0,\xi^j)=\lim_{\mathcal K\in\{\mathcal K\}_\mathcal
A}\,\kappa(x_\mathcal K,\lambda^j_\mathcal K),\label{hr}
\end{equation}
where $j\in I^-$ is arbitrarily chosen.

The case we are thus left with is $I^-\ne\varnothing$. Then,
according to our standing assumptions, $g_{\min}>0$, and therefore
there exists $q\in(0,\infty)$ such that
\begin{equation}
\lambda^j_\mathcal K(\mathbf X)\leqslant q\quad\mbox{for all \ }
\mathcal K\in\{\mathcal K\}_\mathcal A.\label{qq}
\end{equation}
Hence, by (\ref{l11}),
\begin{equation}
\xi^j(\mathbf X)\leqslant q.\label{xi}
\end{equation}

Fix $\varepsilon>0$. Under the assumptions of the lemma, one can
choose a compact neighborhood~$W_{x_0}$ of the point~$x_0$ and a
compact neighborhood~$F$ of the set~$W_{x_0}$ so that
$$W_{x_0}\cap\overline{A_j}=\varnothing,$$
$$F_j:=F\cap\overline{A_j}\ne\varnothing,$$
and
\begin{equation}
\bigl|\kappa(x,y)\bigr|<q^{-1}\varepsilon\quad\mbox{for all \ }
(x,y)\in W_{x_0}\times\complement F.\label{115}
\end{equation}

In the remainder, $\complement_j$ and~$\partial_j$ denote
respectively the complement and the boundary of a set relative
to~$\overline{A_j}$, where $\overline{A_j}$ is treated as a
topological subspace of~$\mathbf X$.

Having observed that $\kappa|_{W_{x_0}\times\overline{A_j}}$\, is
continuous, we proceed to construct a function
$$
\varphi\in\mathbf
C_0(W_{x_0}\times\overline{A_j}\,)
$$
with the following properties:
\begin{equation}
\varphi|_{W_{x_0}\times F_j}=\kappa|_{W_{x_0}\times
F_j},\label{rest}
\end{equation}
\begin{equation}
\bigl|\varphi(x,y)\bigr|\leqslant q^{-1}\varepsilon\quad\mbox{for
all \ } (x,y)\in W_{x_0}\times\complement_jF_j.\label{118}
\end{equation}

To this end, consider a compact neighborhood~$V_j$ of~$F_j$
in~$\overline{A_j}$, and write
$$
f:=\left\{
\begin{array}{cl} \kappa & \mbox{on \ }
W_{x_0}\times\partial_jF_j,\\[2pt]
0 & \mbox{on \ } W_{x_0}\times\partial_jV_j.\\ \end{array} \right.
$$
Note that $E:=(W_{x_0}\times\partial_jF_j)\cup
(W_{x_0}\times\partial_jV_j)$ is a compact subset of the Hausdorff
and compact, hence normal, space $W_{x_0}\times V_j$, and $f$ is
continuous on~$E$. By using the Tietze-Urysohn extension theorem
(see, e.\,g.,~\cite[Th.~0.2.13]{E2}), we deduce from~(\ref{115})
that there exists a continuous function $\hat{f}: \ W_{x_0}\times
V_j\to[-\varepsilon q^{-1},\varepsilon q^{-1}]$ such that
$\hat{f}|_E=f|_E$. Thus, the function in question can be defined as
follows:
$$
\varphi:=\left\{
\begin{array}{cl} \kappa & \mbox{on \ }
W_{x_0}\times F_j,\\[2pt]
\hat{f} & \mbox{on \ } W_{x_0}\times(V_j\setminus F_j),\\[2pt]
0 & \mbox{on \ } W_{x_0}\times\complement_jV_j.
\end{array}
\right.
$$

Furthermore, since the function $\varphi$  is continuous on
$W_{x_0}\times\overline{A_j}$ and has compact support, there exists
a compact neighborhood~$U_{x_0}$ of~$x_0$ such that
\begin{equation}
U_{x_0}\subset W_{x_0}\label{U}
\end{equation}
and
\begin{equation}
\bigl|\varphi(x,y)-\varphi(x_0,y)\bigr|<q^{-1}\varepsilon\quad\mbox{for
all \ } (x,y)\in U_{x_0}\times\overline{A_j}.\label{119}
\end{equation}

Given an arbitrary measure $\nu\in\mathfrak M^+(\,\overline{A_j}\,)$
with the property that $\nu(\mathbf X)\leqslant q$, we conclude
from~\mbox{(\ref{115})\,--\,(\ref{119})} that, for all $x\in
U_{x_0}$,
\begin{equation}
\bigl|\kappa\bigl(x,\nu|_{\complement
F}\bigr)\bigr|\leqslant\varepsilon,\label{121}
\end{equation}
\begin{equation}
\kappa\bigl(x,\nu|_{F}\bigr)=\int\varphi(x,y)\,d\bigl(\nu-\nu|_{\complement
F}\bigr)(y),\label{122}
\end{equation}
\begin{equation}
\Bigl|\int\varphi(x,y)\,d\nu|_{\complement
F}(y)\Bigr|\leqslant\varepsilon,\label{123}
\end{equation}
\begin{equation}
\Bigl|\int
\bigl[\varphi(x,y)-\varphi(x_0,y)\bigr]\,d\nu(y)\Bigr|\leqslant
\varepsilon.\label{124}
\end{equation}

Finally, choose $\mathcal K_0\in\{\mathcal K\}_\mathcal A$ so that,
for all $\mathcal K\succ\mathcal K_0$, there hold $x_\mathcal K\in
U_{x_0}$ and
$$\Bigl|\int\varphi(x_0,y)\,d(\lambda^j_\mathcal
K-\xi^j)(y)\Bigr|<\varepsilon;$$ such a $\mathcal K_0$ exists in
view of~(\ref{l11}) and~(\ref{l112}). Applying now
\mbox{(\ref{121})\,--\,(\ref{124})} to each of~$\lambda_\mathcal
K^j$ and~$\xi^j$, which is possible due to~(\ref{qq})
and~(\ref{xi}), for all $\mathcal K\succ\mathcal K_0$ we therefore
get
\begin{equation*}
\begin{split}
\bigl|\kappa(x_\mathcal K,\lambda_\mathcal
K^j)&-\kappa(x_0,\xi^j)\bigr|\leqslant\bigl|\kappa\bigl(x_\mathcal
K,\lambda_{\mathcal
K}^j\bigl|_{F}\bigr)-\kappa\bigl(x_0,\xi^j\bigl|_{F}\bigr)\bigr|+2\varepsilon\\[7pt]
&{}\leqslant\Bigl|\int\varphi(x_\mathcal K,y)\,d\lambda^j_\mathcal
K(y)-\int\varphi(x_0,y)\,d\xi^j(y)\Bigr|+4\varepsilon\\[7pt]
&{}\leqslant\Bigl|\int\bigl[\varphi(x_\mathcal
K,y)-\varphi(x_0,y)\bigr]\,d\lambda^j_\mathcal K(y)\Bigr|+\Bigl|\int
\varphi(x_0,y)\,d(\lambda_\mathcal
K^j-\xi^j)(y)\Bigr|+4\varepsilon\\[7pt]
&{}\leqslant\varepsilon+\varepsilon+4\varepsilon=6\varepsilon,
\end{split}
\end{equation*}
and (\ref{hr}) follows by letting $\varepsilon$ tend to~$0$. The
proof is complete.\bigskip

{\bf\Large 14. Proof of Theorems 2 and 3}\nopagebreak\medskip

We begin by showing that
\begin{equation} {\rm cap}\,(\mathcal A,a,g)\leqslant\|\hat{\Gamma}(\mathcal
A,a,g)\|^2.  \label{1} \end{equation} To this end,
$\|\hat{\Gamma}(\mathcal A,a,g)\|^2$ can certainly be assumed to be
finite. Then there are $\nu\in\hat{\Gamma}(\mathcal A,a,g)$ and
$\mu\in\mathcal E^0(\mathcal A,a,g)$, the existence of~$\mu$ being
clear from~(\ref{nonzero1}) and Corollary~1.
By~\cite[Lemma~2.3.1]{F1}, the inequality in~(\ref{adm1}) holds
$\mu^i$-almost everywhere. Integrating it with respect to~$\mu^i$
and then summing up over all $i\in I$, in view of $\int
g\,d\mu^i=a_i$ we get
$$ \kappa(\nu,\mu)\geqslant\sum_{i\in I}\,c_i(\nu), $$ hence
$\kappa(\nu,\mu)\geqslant1$ by~(\ref{adm2}), and finally $$
\|\nu\|^2\|\mu\|^2\geqslant1
$$
by the Cauchy-Schwarz inequality. The last relation, being valid for
arbitrary $\nu\in\hat{\Gamma}(\mathcal A,a,g)$ and $\mu\in\mathcal
E^0(\mathcal A,a,g)$, forces (\ref{1}).

The inequality~(\ref{1}) establishes Theorem~2 in the case where
${\rm cap}\,\mathcal A=\infty$.

We are thus left with proving both Theorems~2 and~3 for the case
${\rm cap}\,\mathcal A<\infty$. Then the $\mathcal E(\mathcal
A,a\,{\rm cap}\,\mathcal A,g)$-problem can be considered as well.

Taking (\ref{def}) and (\ref{iden}) into account, we deduce from
Lemmas~5 and~9 with~$a$ replaced by $a\,{\rm cap}\,\mathcal A$ that,
for every $\chi\in\mathcal M'(\mathcal A,a\,{\rm cap}\,\mathcal
A,g)$,
\begin{equation}
\|\chi\|^2={\rm cap}\,\mathcal A
\label{2}
\end{equation}
and there exist unique $\tilde{\eta}_i\in\mathbb R$, $i\in I$, such
that
\begin{equation}
\alpha_ia_i\kappa(x,\chi)\geqslant\tilde{\eta_i}\,g(x)\quad\mbox{n.\,e.~in
\ } A_i,\quad i\in I,\label{tilde1}
\end{equation}
\begin{equation}
\sum_{i\in I}\,\tilde{\eta_i}=1.\label{tilde2}
\end{equation}
Actually,
\begin{equation}
\tilde{\eta_i}=\alpha_i\,{\rm cap}\,\mathcal A^{-1}\,\eta_i(\mathcal
A,a\,{\rm cap}\,\mathcal A,g),\quad i\in I,\label{tilde3}
\end{equation}
where $\eta_i(\mathcal A,a\,{\rm cap}\,\mathcal A,g)$, $i\in I$, are
the numbers uniquely determined in~Sec.~12.

Using the property of subadditivity of $C(\,\cdot\,)$, mentioned
in~Sec.~6.1, and the fact that the potentials of equivalent
in~$\mathcal E$ measures coincide nearly everywhere in~$\mathbf X$,
we conclude from~(\ref{tilde1}) and~(\ref{tilde2}) that
$$
\mathcal M'_\mathcal E(\mathcal A,a\,{\rm cap}\,\mathcal
A,g)\subset\hat{\Gamma}(\mathcal A,a,g).
$$
Together with~(\ref{1}) and~(\ref{2}), this implies that, for every
$\sigma\in\mathcal M'_\mathcal E(\mathcal A,a\,{\rm cap}\,\mathcal
A,g)$,
$${\rm cap}\,\mathcal A=\|\sigma\|^2\geqslant\|\hat{\Gamma}(\mathcal
A,a,g)\|^2\geqslant{\rm cap}\,\mathcal A,$$ which completes the
proof of Theorem~2. The last two relations also yield
$$
\mathcal M'_\mathcal E(\mathcal A,a\,{\rm cap}\,\mathcal
A,g)\subset\hat{\mathcal G}(\mathcal A,a,g).
$$
As both the sides of this inclusion are equivalence classes
in~$\mathcal E$ (see Lemmas~5 and~6), they must actually be equal,
and (\ref{desc}) follows.

Applying Lemma~9$'$ for~$a\,{\rm cap}\,\mathcal A$ in place of~$a$,
we deduce from~(\ref{desc}) that $c_i(\hat{\omega})$, $i\in I$,
satisfying (\ref{adm1}) and~(\ref{adm2}) for
$\nu=\hat{\omega}\in\hat{\mathcal G}(\mathcal A,a,g)$, are
determined uniquely, do not depend on the choice of~$\hat{\omega}$,
and are actually equal to~$\tilde{\eta_i}$. Therefore,
substituting~(\ref{1.3'}) and, subsequently, (\ref{1.3''})
for~$a\,{\rm cap}\,\mathcal A$ in place of~$a$ into~(\ref{tilde3}),
we get~(\ref{snt1}) and (\ref{cnt}). This proves Theorem~3.\bigskip

{\bf\Large 15. Proof of Theorem 5}\nopagebreak\medskip

We start by observing that $\mathcal D(\mathcal A,a,g)$ is nonempty,
contained in an equivalence class in~$\mathcal
E(\,\overline{\mathcal A}\,)$, and satisfies the inclusions
\begin{equation}
\mathcal D(\mathcal A,a,g)\subset\mathcal M(\mathcal A,a\,{\rm
cap}\,\mathcal A,g)\subset\mathcal M'(\mathcal A,a\,{\rm
cap}\,\mathcal A,g)\cap\mathcal E(\,\overline{\mathcal
A},\leqslant\!a\,{\rm cap}\,\mathcal A,g). \label{DD}
\end{equation}
Indeed, this follows from (\ref{D}), Corollary~3, and Lemma~5, the
last two being taken for~$a\,{\rm cap}\,\mathcal A$ in place of~$a$.
Substituting~(\ref{desc2}) into~(\ref{DD}) gives~(\ref{gammain}) as
required.

Since, by (\ref{gammain}), every $\gamma\in\mathcal D(\mathcal
A,a,g)$ is a minimizer in the $\Gamma(\mathcal A,a,g)$-problem, the
claimed relations~(\ref{5.1}) and~(\ref{5.3}) are obtained directly
from Theorem~3 and~4 in view of Definition~5. To show that
$C_i(\mathcal A,a,g)$, $i\in I$, can actually be given by means
of~(\ref{5.2}), one only needs to substitute~$\gamma$ instead
of~$\zeta$ into~(\ref{snt1})~--- which is possible due
to~(\ref{DD})~--- and use Corollary~4.

Assume for a moment that, if $I^-\ne\varnothing$, then $\kappa(x,y)$
is continuous for~$x\ne y$, while $\kappa(\,\cdot\,,y)\to0$ (as
$y\to\infty$) uniformly on compact sets. In order to
establish~(\ref{5.4}), it suffices to apply Lemma~10 (with~$a\,{\rm
cap}\,\mathcal A$ in place of~$a$) to~$\gamma$, which can be done
because of~(\ref{D}), and then substitute~(\ref{5.2}) into the
result obtained.

To prove that $\mathcal D(\mathcal A,a,g)$ is $\mathcal A$-vaguely
compact, fix $(\gamma_s)_{s\in S}\subset\mathcal D(\mathcal A,a,g)$.
Then the inclusion~(\ref{gammain}) and Lemma~7 yield that this net
is $\mathcal A$-vaguely bounded, and hence $\mathcal A$-vaguely
relatively compact. Let $\gamma_0$ denote one of its $\mathcal
A$-vague cluster points, and let $(\gamma_t)_{t\in T}$ be a subnet
of~$(\gamma_s)_{s\in S}$ that converges $\mathcal A$-vaguely
to~$\gamma_0$. In view of~(\ref{D}), the proof will be completed
once we show that
\begin{equation}
\gamma_0\in\mathcal M_0(\mathcal A,a\,{\rm cap}\,\mathcal
A,g).\label{GG}
\end{equation}

By~(\ref{D}), for every $t\in T$ there exist a subnet $(\mathcal
K_{s_t})_{s_t\in S_t}$ of the net $(\mathcal K)_{\mathcal
K\in\{\mathcal K\}_\mathcal A}$ and
$$\lambda_{s_t}\in\mathcal
S(\mathcal K_{s_t},a\,{\rm cap}\,\mathcal A,g),\quad s_t\in S_t,$$
such that $\lambda_{s_t}$ approaches $\gamma_t$ $\mathcal A$-vaguely
as $s_t$ ranges over~$S_t$. Consider the Cartesian product
$\prod\,\{S_t: t\in T\}$~--- that is, the collection of all
functions~$\psi$ on~$T$ with $\psi(t)\in S_t$, and let~$D$ denote
the directed product $T\times\prod\,\{S_t: t\in T\}$ (see,
e.\,g.,~\cite[Chap.~2,~\S~3]{K}). Given $(t,\psi)\in D$, write
$$
\mathcal K_{(t,\psi)}:=\mathcal K_{\psi(t)}\quad\mbox{and}\quad
\lambda_{(t,\psi)}:=\lambda_{\psi(t)}.
$$
Then application of Theorem 4 from \cite[Chap.~2]{K} yields that
$(\lambda_{(t,\psi)})_{(t,\psi)\in D}$ converges $\mathcal
A$-vaguely to~$\gamma_0$. Since, as can be seen from the above
construction, $(\mathcal K_{(t,\psi)})_{(t,\psi)\in D}$ forms a
subnet of $(\mathcal K)_{\mathcal K\in\{\mathcal K\}_\mathcal A}$,
this proves~(\ref{GG}) as required.\bigskip

{\bf\Large 16. Proof of Proposition 2}\nopagebreak\medskip

Consider $\nu\in\mathcal E(\,\overline{\mathcal A}\,)$ and
$\tau_i\in\mathbb R$, $i\in I$, satisfying both the
assumptions~(\ref{desc3'}) and~(\ref{adm2''}), and fix arbitrarily
$\gamma_\mathcal A\in\mathcal D(\mathcal A,a,g)$ and $(\mu_t)_{t\in
T}\in\mathbb M(\mathcal A,a\,{\rm cap}\,\mathcal A,g)$. Since
$\mu^i_t$ is concentrated on~$A_i$ and has finite energy and compact
support, the inequality in~(\ref{desc3'}) holds $\mu^i_t$-almost
everywhere. Integrating it with respect to~$\mu^i_t$ and then
summing up over all $i\in I$, in view of~(\ref{5.1})
and~(\ref{adm2''}) we obtain
$$
2\,\kappa(\mu_t,\nu)\geqslant\|\gamma_\mathcal A\|^2+\|\nu\|^2,\quad
t\in T.
$$
But $(\mu_t)_t\in T$ converges to~$\gamma_\mathcal A$ in the strong
topology of the semimetric space~$\mathcal E(\,\overline{\mathcal
A}\,)$, which is clear from~(\ref{DD}) and Lemma~5 with~$a\,{\rm
cap}\,\mathcal A$ instead of~$a$. Therefore, passing in the
preceding relation to the limit as $t$ ranges over~$T$, we get
$$
\|\nu-\gamma_\mathcal A\|^2=0,
$$
which is a part of the conclusion of the proposition. In turn, the
preceding relation implies that, actually, the right-hand side
in~(\ref{adm2''}) is equal to~$1$, and that $\nu\in\mathcal
M'(\mathcal A,a\,{\rm cap}\,\mathcal A,g)$. Since, in view
of~Theorem~3, the latter means that
$$
R\nu\in\hat{\mathcal G}(\mathcal A,a,g),
$$
the claimed relation (\ref{un}) follows.\bigskip

{\bf\Large 17. Proof of Theorem 6}\nopagebreak\medskip

To establish~(\ref{6.1}), fix $\mu\in\mathcal E(\mathcal A,a,g)$.
Under the assumptions of the theorem, either $g_{\min}>0$, and
consequently $\mu^i(\mathbf X)<\infty$ for all $i\in I$, or $\mathbf
X$ is countable at infinity; in any case, every~$A_i$, $i\in I$, is
contained in a countable union of $\mu^i$-integrable sets.
Therefore, by~\cite{B2,E2} (cf.~the appendix below),
\begin{align*}
\int g\,d\mu^i&=\lim_{n\in\mathbb N}\,\int
g\,d\mu^{i}_{\mathcal A_n},&&\!\!\!\!\!\!\!\!\!\!\!\!\!\!\!\!\!\!\!\!\!\!\!\!\!\!\!\!\!\!i\in I,\\
\kappa(\mu^i,\mu^j)&=\lim_{n\in\mathbb N}\,\kappa(\mu^{i}_{\mathcal
A_n},\mu^{j}_{\mathcal
A_n}),&&\!\!\!\!\!\!\!\!\!\!\!\!\!\!\!\!\!\!\!\!\!\!\!\!\!\!\!\!\!\!i,\,j\in
I,
\end{align*}
where $\mu^{i}_{\mathcal A_n}$ denotes the trace of $\mu^{i}$
upon~$A^i_n$. Now, applying the same arguments as in the proof of
Lemma~2, but with the preceding two relations instead of~(\ref{mon})
and~(\ref{mon'}), we arrive at~(\ref{6.1}) as required.

By~(\ref{nonzero1}) and~(\ref{6.1}), for every $n\in\mathbb N$,
${\rm cap}\,(\mathcal A_n,a,g)$ can certainly be assumed to be
nonzero. Suppose moreover that ${\rm cap}\,(\mathcal A,a,g)$ is
finite; then, by~(\ref{increas}), so is ${\rm cap}\,(\mathcal
A_n,a,g)$. Hence, according to Theorem~5, there exists
\begin{equation}
\gamma_n:=\gamma_{\mathcal A_n}\in\mathcal D(\mathcal
A_n,a,g).\label{17.1}\end{equation}

Observe that $R\gamma_n$ is a minimizer in the
$\hat{\Gamma}(\mathcal A_n,a,g)$-problem, which is clear
from~(\ref{desc}), (\ref{desc2}), and~(\ref{gammain}). Since,
furthemore,
$$\hat{\Gamma}(\mathcal A_{n+1},a,g)\subset\hat{\Gamma}(\mathcal A_n,a,g),$$
application of Lemma~1 to $\mathcal H=\hat{\Gamma}(\mathcal
A_n,a,g)$, $\nu=R\gamma_{n+1}$, and $\lambda=R\gamma_n$ gives
$$
\|\gamma_{n+1}-\gamma_n\|^2\leqslant\|\gamma_{n+1}\|^2-\|\gamma_n\|^2.
$$
Also note that $\|\gamma_n\|^2$, $n\in\mathbb N$, is a Cauchy
sequence in~$\mathbb R$, because, by~(\ref{6.1}), its limit exists
and, being equal to~${\rm cap}\,\mathcal A$, is finite. The
preceding inequality therefore yields that $(\gamma_n)_{n\in\mathbb
N}$ is a strong Cauchy sequence in the semimetric space~$\mathcal
E(\,\overline{\mathcal A}\,)$.

Besides, since ${\rm cap}\,\mathcal A_n\leqslant{\rm cap}\,\mathcal
A$, we derive from~(\ref{gammain}) that
$$(\gamma_n)_{n\in\mathbb N}\subset\mathcal E(\,\overline{\mathcal
A}\,,\leqslant\!a\,{\rm cap}\,\mathcal A,g).$$ Hence, by Theorem~7,
there exists an $\mathcal A$-vague cluster point~$\gamma$
of~$(\gamma_n)_{n\in\mathbb N}$, and
$$
\lim_{n\in\mathbb N}\,\|\gamma_n-\gamma\|^2=0.
$$
Let $(\gamma_t)_{t\in T}$ denote a subnet of the sequence
$(\gamma_n)_{n\in\mathbb N}$ that converges $\mathcal A$-vaguely and
strongly to~$\gamma$. We next proceed to show that
\begin{equation}
\gamma\in\mathcal D(\mathcal A,a,g). \label{theta}\end{equation}

For every $t\in T$, consider the ordered family $\{\mathcal
K_t\}_{\mathcal A_t}$ of all compact condensers $\mathcal
K_t\prec\mathcal A_t$. By~(\ref{17.1}), there exist a subnet
$(\mathcal K_{s_t})_{s_t\in S_t}$ of $(\mathcal K_t)_{\mathcal
K_t\in\{\mathcal K_t\}_{\mathcal A_t}}$ and
$$\lambda_{s_t}\in\mathcal S(\mathcal K_{s_t},a\,{\rm
cap}\,\mathcal K_{s_t},g)$$ such that $(\lambda_{s_t})_{s_t\in S_t}$
converges both strongly and $\mathcal A$-vaguely to~$\gamma_t$.
Consider the Cartesian product $\prod\,\{S_t: t\in T\}$, that is,
the collection of all functions~$\psi$ on~$T$ with $\psi(t)\in S_t$,
and let $D$ denote the directed product $T\times\prod\,\{S_t: t\in
T\}$. Given $(t,\psi)\in D$, write
$$
\mathcal K_{(t,\psi)}:=\mathcal K_{\psi(t)}\quad\mbox{and}\quad
\lambda_{(t,\psi)}:=\lambda_{\psi(t)}.
$$
Then application of Theorem 4 from~\cite[Chap.~2]{K} yields that
$(\lambda_{(t,\psi)})_{(t,\psi)\in D}$ converges both strongly and
$\mathcal A$-vaguely to~$\gamma$. Since $(\mathcal
K_{(t,\psi)})_{(t,\psi)\in D}$ is easily checked to form a subnet
of~$(\mathcal K)_{\mathcal K\in\{\mathcal K\}_\mathcal A}$, this
proves~(\ref{theta}) as required.

What is finally left is to prove~(\ref{6.3}). By Corollary~9, for
every $n\in\mathbb N$ one can choose a compact condenser $\mathcal
K^0_n\prec\mathcal A_n$ so that
$$
\bigl|C_i(\mathcal A_n,a,g)-C_i(\mathcal
K^0_n,a,g)\bigr|<n^{-1},\quad i\in I.
$$
This $\mathcal K^0_n$ can certainly be chosen so large that the
sequence obtained, $(\mathcal K^0_n)_{n\in\mathbb N}$, forms a
subnet of $(\mathcal K)_{\mathcal K\in\{\mathcal K\}_\mathcal A}$;
therefore, repeated application of Corollary~9 yields
$$\lim_{n\in\mathbb N}\,C_i(\mathcal
K^0_n,a,g)=C_i(\mathcal A,a,g).$$ This leads to (\ref{6.3}) when
combined with the preceding relation.\bigskip

{\bf\Large 18. Acknowledgments}\nopagebreak\medskip

The author is greatly indebted to Professors W.~Hansen, E.~Saff, and
W.~Wendland for several helpful comments concerning this study, and
to Professor~B.~Fuglede for drawing the author's attention to the
articles~\cite{F2} and~\cite{E1}.\bigskip

{\bf\Large 19. Appendix}\nopagebreak\medskip

Let $\nu\in\mathfrak M^+(\mathbf X)$ be given. As in~\cite[Chap.~4,
\S~4.7]{E2}, a set~$E\subset\mathbf X$ is called
\mbox{$\nu$-$\sigma$}-{\it fi\-ni\-te\/} if it can be written as a
countable union of $\nu$-integr\-able sets.

The following assertion, related to the theory of measures and
integration, has been used in Sec.~17. Although it is not difficult
to deduce it from~\cite{B2,E2}, we could not find there a proper
reference.\medskip

{\bf Lemma 11.} {\it Consider a lower semicontinuous function~$\psi$
on~$\mathbf X$ such that $\psi\geqslant0$ unless the space~$\mathbf
X$ is compact, and let $E$ be the union of an increasing sequence of
$\nu$-measur\-able sets~$E_n$, $n\in\mathbb N$. If moreover $E$ is
\mbox{$\nu$-$\sigma$}-fi\-ni\-te, then}
$$
\int\psi\,d\nu_E=\lim_{n\in\mathbb N}\,\int\psi\,d\nu_{E_n}.
$$

{\bf Proof.} We can certainly assume $\psi$ to be nonnegative, for
if not, we replace $\psi$ by a function~$\psi'$ obtained by adding
to~$\psi$ a suitable constant~$c>0$:
$$
\psi'(x):=\psi(x)+c\geqslant0,
$$
which is always possible since a lower semicontinuous function is
bounded from below on a compact space. Then, for every
$\nu$-measurable and \mbox{$\nu$-$\sigma$}-fi\-ni\-te set~$Q$,
\begin{equation}
\int\psi\,d\nu_Q=\int\psi\varphi_Q\,d\nu, \label{Q}
\end{equation}
where $\varphi_Q(x)$ equals~$1$ if $x\in Q$, and~$0$ otherwise.
Indeed, this can be concluded from~\cite[Chap.~4, \S~4.14]{E2} (see
Propositions~4.14.1 and~4.14.6).

On the other hand, since $\psi\varphi_{E_n}$, $n\in\mathbb N$, are
nonnegative and form an increasing sequence with the upper
envelope~$\psi\varphi_E$, \cite[Prop.~4.5.1]{E2} gives
$$
\int\psi\varphi_E\,d\nu=\lim_{n\in\mathbb
N}\,\int\psi\varphi_{E_n}\,d\nu.
$$
Applying (\ref{Q}) to both the sides of this equality, we obtain the
lemma.
\bigskip
\flushleft

{\small Institute of Mathematics\\
National Academy of Sciences of Ukraine\\
3 Tereshchenkivska Str.\\
01601, Kyiv-4, Ukraine\\
e-mail: natalia.zorii@gmail.com}

\end{document}